\PassOptionsToPackage{unicode}{hyperref}
\PassOptionsToPackage{hyphens}{url}
\PassOptionsToPackage{dvipsnames,svgnames,x11names}{xcolor}
\documentclass[
  11pt]{article}

\usepackage{amsmath,amssymb}
\usepackage{iftex}
\ifPDFTeX
  \usepackage[T1]{fontenc}
  \usepackage[utf8]{inputenc}
  \usepackage{textcomp} 
\else 
  \usepackage{unicode-math}
  \defaultfontfeatures{Scale=MatchLowercase}
  \defaultfontfeatures[\rmfamily]{Ligatures=TeX,Scale=1}
\fi
\usepackage{lmodern}
\ifPDFTeX\else  
\fi
\IfFileExists{upquote.sty}{\usepackage{upquote}}{}
\IfFileExists{microtype.sty}{
  \usepackage[]{microtype}
  \UseMicrotypeSet[protrusion]{basicmath} 
}{}
\makeatletter
\@ifundefined{KOMAClassName}{
  \IfFileExists{parskip.sty}{%
    \usepackage{parskip}
  }{
    \setlength{\parindent}{0pt}
    \setlength{\parskip}{6pt plus 2pt minus 1pt}}
}{
  \KOMAoptions{parskip=half}}
\makeatother
\usepackage{xcolor}
\makeatletter
\ifx\paragraph\undefined\else
  \let\oldparagraph\paragraph
  \renewcommand{\paragraph}{
    \@ifstar
      \xxxParagraphStar
      \xxxParagraphNoStar
  }
  \newcommand{\xxxParagraphStar}[1]{\oldparagraph*{#1}\mbox{}}
  \newcommand{\xxxParagraphNoStar}[1]{\oldparagraph{#1}\mbox{}}
\fi
\ifx\subparagraph\undefined\else
  \let\oldsubparagraph\subparagraph
  \renewcommand{\subparagraph}{
    \@ifstar
      \xxxSubParagraphStar
      \xxxSubParagraphNoStar
  }
  \newcommand{\xxxSubParagraphStar}[1]{\oldsubparagraph*{#1}\mbox{}}
  \newcommand{\xxxSubParagraphNoStar}[1]{\oldsubparagraph{#1}\mbox{}}
\fi
\makeatother

\usepackage{longtable,booktabs,array}
\usepackage{calc} 
\usepackage{etoolbox}
\makeatletter
\patchcmd\longtable{\par}{\if@noskipsec\mbox{}\fi\par}{}{}
\makeatother
\IfFileExists{footnotehyper.sty}{\usepackage{footnotehyper}}{\usepackage{footnote}}
\makesavenoteenv{longtable}
\usepackage{graphicx}
\makeatletter
\def\maxwidth{\ifdim\Gin@nat@width>\linewidth\linewidth\else\Gin@nat@width\fi}
\def\maxheight{\ifdim\Gin@nat@height>\textheight\textheight\else\Gin@nat@height\fi}
\makeatother
\setkeys{Gin}{width=\maxwidth,height=\maxheight,keepaspectratio}
\makeatletter
\def\fps@figure{htbp}
\makeatother

\makeatletter
\@ifpackageloaded{caption}{}{\usepackage{caption}}
\AtBeginDocument{%
\ifdefined\contentsname
  \renewcommand*\contentsname{Table of contents}
\else
  \newcommand\contentsname{Table of contents}
\fi
\ifdefined\listfigurename
  \renewcommand*\listfigurename{List of Figures}
\else
  \newcommand\listfigurename{List of Figures}
\fi
\ifdefined\listtablename
  \renewcommand*\listtablename{List of Tables}
\else
  \newcommand\listtablename{List of Tables}
\fi
\ifdefined\figurename
  \renewcommand*\figurename{Figure}
\else
  \newcommand\figurename{Figure}
\fi
\ifdefined\tablename
  \renewcommand*\tablename{Table}
\else
  \newcommand\tablename{Table}
\fi
}
\@ifpackageloaded{float}{}{\usepackage{float}}
\floatstyle{ruled}
\@ifundefined{c@chapter}{\newfloat{codelisting}{h}{lop}}{\newfloat{codelisting}{h}{lop}[chapter]}
\floatname{codelisting}{Listing}

\makeatother
\makeatletter
\@ifpackageloaded{caption}{}{\usepackage{caption}}
\@ifpackageloaded{subcaption}{}{\usepackage{subcaption}}
\makeatother

\ifLuaTeX
  \usepackage{selnolig}  
\fi
\usepackage[]{natbib}
\usepackage{setspace}
\usepackage{bookmark}

\IfFileExists{xurl.sty}{\usepackage{xurl}}{} 
\hypersetup{
  pdftitle={Specification LASSO and a Flexible Characteristics-Based Asset Pricing Model},
  pdfauthor={Chaohua Dong; Shuyi Ge; Shaoran Li; Wen Su},
  pdfkeywords={Interaction of covariates; LASSO; Model selection; Variable selection},
  colorlinks=true,
  linkcolor={blue},
  filecolor={Maroon},
  citecolor={Blue},
  urlcolor={Blue},
  pdfcreator={LaTeX}}

\usepackage{colortbl}
\usepackage{multirow}
\usepackage{bm}
\usepackage{amsthm}
\theoremstyle{plain}
\newtheorem{thm}{Theorem}[section]

\newtheorem{lem}[thm]{Lemma}

\theoremstyle{definition}
\newtheorem{defn}{Definition}[section]

\newtheorem{exmp}{Example}[section]
\theoremstyle{remark}

\theoremstyle{remark}
\newtheorem{assu}{Assumption}
\newcommand{\rmk}{\noindent{\textit{\textbf{Remark}}}. }
\newcommand{\prf}{\noindent{\textit{\textbf{Proof}}}. }
\usepackage{longtable}

\usepackage{fancyhdr}

\usepackage{chapterbib}
\usepackage{etoc}
\usepackage{mathtools,nccmath}
\usepackage{adjustbox}

\etocsettagdepth{main}{none}
\etocsettagdepth{appendix}{subsubsection}

\newsavebox{\SLASSOTableBox}
\BeforeBeginEnvironment{tabular}{\begin{lrbox}{\SLASSOTableBox}}
\AfterEndEnvironment{tabular}{%
  \end{lrbox}%
  \ifdim\wd\SLASSOTableBox>\linewidth
    \resizebox{\linewidth}{!}{\usebox{\SLASSOTableBox}}%
  \else
    \usebox{\SLASSOTableBox}%
  \fi}

\newlength{\SLASSOColumnWidth}
\let\SLASSOOriginalLongtable\longtable
\RenewDocumentCommand{\longtable}{o m}{%
  \ifstrequal{#2}{>{\itshape}p{2.2cm}p{5.5cm}>{\itshape}p{2.2cm}p{5.5cm}}{%
    \setlength{\SLASSOColumnWidth}{\dimexpr\linewidth-8\tabcolsep\relax}%
    \setstretch{.93}%
    \setlength{\LTpre}{12pt}%
    \setlength{\LTpost}{6pt}%
    \IfValueTF{#1}{%
      \SLASSOOriginalLongtable[#1]{>{\itshape}p{.15\SLASSOColumnWidth}p{.35\SLASSOColumnWidth}>{\itshape}p{.15\SLASSOColumnWidth}p{.35\SLASSOColumnWidth}}}{%
      \SLASSOOriginalLongtable{>{\itshape}p{.15\SLASSOColumnWidth}p{.35\SLASSOColumnWidth}>{\itshape}p{.15\SLASSOColumnWidth}p{.35\SLASSOColumnWidth}}}%
  }{%
    \IfValueTF{#1}{\SLASSOOriginalLongtable[#1]{#2}}{\SLASSOOriginalLongtable{#2}}%
  }%
}

\newif\ifSLASSOAppendix
\AtBeginEnvironment{equation}{%
  \ifSLASSOAppendix
    \ifnum\value{equation}=70\relax\small\fi
  \fi}

\begin{document}
\def\spacingset#1{\renewcommand{\baselinestretch}{#1}\small\normalsize}
\spacingset{1.5}
\etocdepthtag{main}
\title{\bfseries Specification LASSO and a Flexible Characteristics-Based Asset Pricing Model\footnote{\scriptsize Forthcoming at {\it Journal of Business \& Economic Statistics}.}}
\author{%
  \parbox{\dimexpr\textwidth-2\tabcolsep\relax}{\centering
    Chaohua Dong$^1$, Shuyi Ge$^2$, Shaoran Li$^3$, Wen Su$^4$\footnote{\scriptsize Emails: dchaohua@zuel.edu.cn, sg751$\_$shuyige@outlook.com, lishaoran@pku.edu.cn, wen.su@maths.ox.ac.uk}\\[0.7em]
    \small\itshape
    \textsuperscript{1}School of Statistics and Mathematics,
    Zhongnan University of Economics and Law\\
    \textsuperscript{2}School of Finance, Nankai University\\
    \textsuperscript{3}School of Economics, Peking University\\
    \textsuperscript{4}Mathematical Institute, University of Oxford}}
\date{}
\maketitle

\textbf{Abstract}: 
This paper studies a partially linear semiparametric additive model that allows for interactions among variables in high-dimensional settings.
We propose the Specification LASSO (S-LASSO), a two-step method that integrates LASSO and Adaptive Group LASSO to simultaneously perform variable selection and model specification.
Theoretically, we establish that S-LASSO possesses oracle properties in terms of selection consistency and asymptotic normality.
Monte Carlo simulations confirm that S-LASSO outperforms existing alternatives across a variety of scenarios.
Empirically, we apply S-LASSO to a characteristics-based asset pricing model, exploring the linear and nonlinear effects of asset-specific characteristics on cross-sectional returns, as well as their interactions with firm size.


\section{Introduction}
\label{sec1}

In the data-rich era, researchers increasingly face two fundamental challenges: \emph{variable selection} and \emph{model specification}. 
The \emph{variable selection} problem arises from the abundance of available covariates, especially in high-dimensional settings where the number of predictors \( P \) may exceed the number of observations \( n \). 
At the same time, the \emph{model specification} challenge concerns the choice of an appropriate functional form from an expanding ``model zoo.'' 
Parametric models are widely used for their simplicity and interpretability, but they are vulnerable to misspecification.
Nonparametric models offer greater flexibility by relaxing functional assumptions, yet they often suffer from the \emph{curse of dimensionality}, which limits their applicability when the number of covariates is large.
Modern machine learning and deep learning methods provide powerful predictive tools, but they typically lack explicit model representations and therefore offer limited interpretability regarding the underlying economic or structural relationships.

A central question in asset pricing is which firm characteristics contain independent information about returns, and in what functional forms they do so.
This question is challenging because the literature has documented a ``zoo'' of characteristics that appear to predict the cross-section of returns, creating both variable-selection and model-specification problems.

In this paper, we study a high-dimensional partially linear semiparametric additive model with interactions, building on LASSO and adaptive group LASSO \citep{Tibshirani1996, Zou2006, Yuan2006, Wang2008}. 
We propose a new specification and estimation framework, \emph{Specification LASSO} (S-LASSO). 
Under mild regularity conditions, S-LASSO consistently selects relevant variables and identifies whether they enter linearly, nonlinearly, or through interactions.
It therefore addresses variable selection and model specification jointly in high-dimensional environments.
In our empirical application, S-LASSO yields a more efficiently specified characteristics-based asset pricing model with better out-of-sample performance.

The method begins with a comprehensive set of candidate components, including linear terms, nonlinear terms, and pairwise interactions. 
Nonlinear effects are approximated by orthogonal series expansions, interactions by tensor products of basis functions, and the resulting components are selected and estimated through a two-step framework that combines LASSO and adaptive group LASSO.

Our framework is built on a sparse high-dimensional regression structure. 
Suppose we observe $\left(\bm{Y}, \bm{P}\right)$, where $\bm{Y}$ is an $n \times 1$ vector of dependent variables and $\bm{P}$ is an $n \times P$ matrix of potential covariates $\left(\bm{X}_{1}, \bm{X}_{2}, \dots, \bm{X}_{P}\right)$, allowing $P > n$.
We assume that only a subset $\bm{Q}$ of size $Q$ is relevant, where $Q$ is fixed while $P$ diverges as $n \to \infty$. 
For notational convenience, we assume that the first $Q$ variables in $\bm{P}$ are relevant. 
Specifically,
\begin{equation}
E\left(Y_i \mid X_{1i}, \dots, X_{Pi}\right) = \theta + h\left(X_{1i}, \dots, X_{Qi}\right),
\end{equation}
where $\theta$ is a constant and $h(\cdot)$ is an unknown multivariate function. 
The idiosyncratic error is $\epsilon_i = Y_i - E\left(Y_i \mid X_{1i}, \dots, X_{Pi}\right)$.

A common starting point is the additive nonparametric model
$
h(\bm{Q}_{i})=\sum_{q=1}^{Q}f_{q}(X_{qi}),
$
where \(f_q(\cdot)\) is an unknown univariate function. 
Such models have been widely studied by \cite{hastie1990exploring}, \cite{linton1997miscellanea}, \cite{linton2000efficient}, and \cite{linton1996estimation}. 
Although the additive structure alleviates the curse of dimensionality, it can be inefficient when some covariates have linear effects, since nonparametric estimators generally converge more slowly than $O(n^{-1/2})$.

This motivates the partially linear additive semiparametric model
\begin{equation}
h(\bm{Q_{i}})=\sum_{l=1}^{L}\beta_{l}X_{il}+\sum_{q=L+1}^{Q}f_{q}(X_{iq}),
\end{equation}
which separates linear and nonlinear effects. 
The linear coefficients can then be estimated at the rate $O(n^{-1/2})$; see \citet{wang2007spline} and \cite{ma2011spline}. 
Related models are also studied by \citet{li2000efficient}, \citet{fan2003kernel}, and \citet{liang2008additive}.

However, these additive models do not account for interactions among covariates.
As the following examples show, pairwise interactions are common in economics and finance, so ignoring them can be restrictive.

\begin{exmp}
    In macroeconomics, many production functions incorporate interactions between inputs. For example, the generalized Cobb-Douglas production function can be written as
    \(
    \log Y = c + \sum_{\alpha=1}^d \sum_{\beta=1}^d c_{\alpha \beta} \log \left(\frac{X_\alpha + X_\beta}{2}\right),
    \)
    where \(X_\alpha\) and \(X_\beta\) denote production inputs and \(c_{\alpha \beta}\) captures their interaction effect.
\end{exmp}

\begin{exmp}
    In microeconomics, \citet{deaton1980economics} document a utility model for a household $Y$ that includes interactions between variables such as eating and drinking ($X_E, X_D$), housing and fuel ($X_H, X_F$), and television and sports ($X_T, X_S$):
    \(
    Y = m_{ED}(X_E, X_D) + m_{HF}(X_H, X_F) + m_{TS}(X_T, X_S) + \epsilon.
    \)
\end{exmp}

\begin{exmp}
    In environmental studies, \citet{Dong2021JoE} examine the effects of $\mathrm{CO}_2$ and solar irradiance on global sea level rise and specify
    \(
    Y_{\mathrm{SL}} = m(\tau, X_{\mathrm{CO}_2}, X_{\mathrm{SI}}) + \epsilon,
    \)
    documenting interaction effects among time, $\mathrm{CO}_2$, and solar irradiance.
\end{exmp}

\begin{exmp}
    In financial studies, \citet{ALL2021RFS} estimate the implied volatility surface using
    \(
    \sigma(T, K) = \sum_{j=0}^J \sum_{i=0}^{L_j} \beta_{ij} T^i K^j,
    \)
    where $K$ is log-moneyness and $T$ is time to maturity.
\end{exmp}

Motivated by these examples, we model the unknown multivariate function as a combination of linear terms, nonlinear terms, and pairwise interactions.
\citet{sperlich2002nonparametric} emphasize the importance of interactions in additive models and propose marginal integration methods for estimation and testing. 
Their approach, however, is computationally demanding and restricted to low-dimensional settings.

S-LASSO fills this gap by modeling and selecting pairwise interactions in high-dimensional settings. 
Pairwise interactions offer a useful compromise between interpretability and flexibility, while the framework can be extended to higher-order interactions if needed. 
We establish oracle properties for S-LASSO: it consistently selects relevant components and their functional forms, and the estimators of selected components have the same asymptotic distribution as their ordinary least squares counterparts. 
Monte Carlo experiments confirm S-LASSO's strong performance even when \( P > n \).

Our empirical application studies \emph{which} firm characteristics explain the cross-section of asset returns, and \emph{how} they do so. 
We extend the dataset of \citet{GHZ2017RFS} and collect 102 firm-level characteristics. 
Building on this literature, \citet{freyberger2017dissecting} approximate the nonlinear effects of 62 characteristics with B-splines, apply adaptive group LASSO for selection, and introduce interaction terms by multiplying firm size with the remaining characteristics.

We depart from \citet{freyberger2017dissecting} in two ways.
First, we allow linear, nonlinear, and interactive effects to coexist.
Second, we allow interactions with firm size to take general bivariate forms rather than simple multiplicative ones.
Empirically, S-LASSO outperforms the benchmark of \citet{freyberger2017dissecting} both in sample and out of sample, reflecting its ability to separate linear from nonlinear effects and accommodate more flexible size interactions.

\textbf{Literature Review}

Variable selection in high-dimensional settings has been studied more extensively than model specification for unknown multivariate functions. 
One of the most widely used approaches is LASSO (Least Absolute Shrinkage and Selection Operator). 
\citet{Tibshirani1996} introduced LASSO to perform regularization and variable selection in linear models:
$
\min_{\bm{\alpha}} \sum_{i=1}^{n}\left(Y_i-\sum_{j=1}^{W}\alpha _jX_{ij}\right)^2+\lambda _n\sum_{j=1}^{W}|\alpha _j| ,
$
where $\lambda_{n}$ is a tuning parameter. 
LASSO can estimate coefficients and set irrelevant ones to zero even when $P \gg n$.
A necessary condition for sign-consistent selection is the irrepresentable condition; see \citet{zhao2006model} and \citet{Zou2006}. 
We revisit this condition in Section \ref{sec4}.

To obtain $\sqrt{n}$-consistency, \citet{Zou2006} proposed the Adaptive LASSO:
\begin{equation*}
\min_{\bm{\beta}} \sum_{i=1}^{n}\left(Y_i-\sum_{j=1}^{W}\beta _jX_{ij}\right)^2+\lambda _n\sum_{j=1}^{W}\widehat{w}_j|\beta _j| ,
\end{equation*}
where $\widehat{w}_j$ is typically chosen as $\widehat{w}_j = |\widehat{\alpha}_j|^{-\gamma}$ for some $\gamma > 0$, with $\widehat{\alpha}_j$ a preliminary consistent estimator.

Traditional LASSO methods mainly focus on univariate effects. 
However, \cite{FKLL2016} argue that interactions are also important in many applications. 
They study
\(
Y_i = \beta_0 + \sum_{j=1}^{d} \beta_j X_i^{(j)} + \sum_{j<k} \beta_{jk}X_i^{(j)}X_i^{(k)} + \varepsilon_i,
\)
and establish sure screening and oracle inequalities for both univariate effects and interactions.

Another related method is \emph{COSSO}, proposed by \citet{LZ2006AoS} under the smoothing spline ANOVA framework:
\begin{equation}
h(\bm{X}_i) = b + \sum_{j=1}^{d} f_j(X_i^{(j)}) + \sum_{j<k} f_{jk}(X_i^{(j)}, X_i^{(k)}) + \cdots.
\end{equation}
COSSO is flexible, but it assumes a fixed dimension, does not distinguish linear and nonlinear components, and does not yield an explicit parameterized representation. 
By contrast, our method accommodates linear terms, performs variable and model selection, and is applicable in high-dimensional settings.

\citet{Huang2010} propose a selection and estimation method for additive nonparametric models based on the group LASSO and adaptive group LASSO. 
They approximate each unknown function by B-spline bases,
$f_{nj}(x) = \sum_{k=1}^{m_n} \beta_{jk} \phi_k(x)$,
and then minimize
\begin{equation}
L_n(\mu, \beta_n) = \sum_{i=1}^n \left(Y_i - \mu - \sum_{j=1}^P \sum_{k=1}^{m_n} \beta_{jk} \phi_k(X_{ij})\right)^2 + \lambda_n \sum_{j=1}^P \widehat{w}_{nj} \|\bm{\beta}_{nj}\|_2,
\end{equation}
where $\lambda_n$ is a tuning parameter, $\|\bm{\beta}_{nj}\|_2$ is the Euclidean norm of $\bm{\beta}_{nj} = (\beta_{j1}, \dots, \beta_{jm_n})^\intercal$, and $\widehat{w}_{nj}$ is constructed from an initial estimator. 
They show that adaptive group LASSO can outperform COSSO in selecting relevant variables in semiparametric additive models.

The contribution of this paper is threefold. 
First, we propose a more general and flexible modeling framework that incorporates linear components with faster convergence rates, nonparametric components to capture nonlinear effects, and pairwise interaction terms to allow for richer functional relationships. 
Second, we develop the S-LASSO method, which simultaneously performs variable selection, model specification, and estimation.
We establish its oracle properties in terms of variable selection consistency and show the asymptotic normality of the estimators. 
Finally, we empirically demonstrate the presence and importance of interactive effects between firm size and other characteristics in characteristics-based asset pricing models.

The rest of the paper is organized as follows. 
Section \ref{sec2} introduces the model framework.
Section \ref{sec3} describes the S-LASSO procedure. 
Section \ref{sec4} establishes the theoretical results. 
Section \ref{sec5} reports simulation evidence, and Section \ref{sec6} presents the empirical study. 
Section \ref{sec7} concludes. 
All proofs and additional materials are collected in the Appendices.
\section{Model Setup}
\label{sec2}

\subsection{Partially Linear Additive Semiparametric Model}

Suppose that we observe a sample $\left(\bm{Y}, \bm{P}\right)$ and study the following sparse model:
\begin{equation}
Y_i = \theta + h\left(X_{1i}, \dots, X_{Qi}\right) + \epsilon_i,
\end{equation}
where $\theta$ is a constant, $\epsilon_i$ is an idiosyncratic error term, and $h(\cdot)$ is an unknown multivariate function. 
In vector form, the model can be written as $\bm{Y} = \bm{\theta} + \bm{h}(\bm{Q}) + \bm{U}$, where $\bm{\theta} = \theta \bm{1}$, $\bm{U} = (\epsilon_1, \dots, \epsilon_n)^\intercal$ is the error vector, and
$
\bm{h}(\bm{Q}) = \left(h\left(X_{11}, \dots, X_{Q1}\right), \dots, h\left(X_{1n}, \dots, X_{Qn}\right)\right)^\intercal
$.

In this paper, we specify \( h(\cdot) \) as a partially linear additive semiparametric model with interaction terms:
\begin{equation}\label{equation 4}
\bm{h}(\bm{Q}) = 
\underbrace{\sum_{l \in \boldsymbol{L}} \alpha_l \boldsymbol{X}_l}_{\text{linear}} + 
\underbrace{\sum_{r \in \boldsymbol{R}} m_r\left(\boldsymbol{X}_r\right)}_{\text{nonlinear}} + 
\underbrace{\sum_{(s,s') \in \boldsymbol{S}} m_{ss'}\left(\boldsymbol{X}_s, \boldsymbol{X}_{s'}\right)}_{\text{interactive}},
\end{equation}
where \(\bm{Q}\) is organized into component sets corresponding to linear effects, nonlinear effects, and interaction effects, and all of which are inferred from the data.
Specifically, we define
\begin{equation}\label{eqn:sets}
\begin{aligned}
    &\boldsymbol{L} = \left\{ l : \boldsymbol{X}_l  \text{ has a linear effect on } \bm{Y} \right\}, \\
    &\boldsymbol{R} = \left\{ r : \boldsymbol{X}_r \text{ has a nonlinear effect on } \bm{Y} \right\}, \\
    &\boldsymbol{S} = \left\{ (s, s') : s < s', \boldsymbol{X}_s,\boldsymbol{X}_{s'} \text{ have an interactive effect on } \bm{Y} \right\}.
\end{aligned}
\end{equation}
The ordering condition \(s < s'\) avoids double-counting interaction effects.
 
In this setup, $m_{r}(\bm{X}_r)$ denotes an unknown univariate function, and $m_{ss'}(\bm{X}_s, \bm{X}_{s'})$ denotes an unknown bivariate function. 
The cardinalities of the sets in \autoref{eqn:sets} are $|\bm{L}| = L$, $|\bm{R}| = R$, and $|\bm{S}| = S$, where $S$ is the number of distinct variable pairs. 
Each set is unknown. 
We also define the complement sets of irrelevant variables, $\bm{L}^C$, $\bm{R}^C$ and $\bm{S}^C$.

The decomposition in \autoref{equation 4} mitigates the curse of dimensionality by restricting nonparametric estimation to univariate and bivariate functions while retaining substantial flexibility.
Relative to conventional additive models with only univariate components, and relative to \citet{freyberger2017dissecting}, which restrict interactions to multiplicative forms involving firm size, our model allows general interactions among covariates and is therefore more flexible.
The inclusion of linear components also improves estimation efficiency and reduces computational burden when linear approximations are adequate.
Although we focus on pairwise interactions for interpretability and tractability, the framework can be extended to higher-order interactions if needed.

\subsection{Parameterization}

To formalize the method and establish the theoretical results, we introduce notation and define the parameterization. 
Following \cite{Dong2021JoE}, we assume that the \(j^{\text{th}}\) potential covariate \(\bm{X}_j\) is an \(n \times 1\) vector
\(
\bm{X}_{j} = \left(X_{j1}, X_{j2}, \dots, X_{jn}\right)^{\intercal}\) for $j = 1, 2, \dots, P$, which takes values in \(\Omega^n\), where \(\Omega = \mathbb R\) in this paper.\footnote{For theoretical convenience, $\Omega$ is sometimes restricted to be a compact subset of $\mathbb R$, especially for sieve methods. For example, see \cite{AiandChen2003ECTA}.}

The unknown functions \(m_r : \Omega \to \mathbb{R}\) and \(m_{ab} : \Omega^2 \to \mathbb{R}\) are extended to their multivariate forms as
\begin{equation}
\begin{aligned}
    & \bm{m}_r\left(\bm{X}_{r}\right) = \left(m_r\left(X_{r1}\right), m_r\left(X_{r2}\right), \dots, m_r\left(X_{rn}\right)\right)^{\intercal}, \\
    & \bm{m}_{ab}\left(\bm{X}_{a}, \bm{X}_{b}\right) = \left(m_{ab}\left(X_{a1}, X_{b1}\right), m_{ab}\left(X_{a2}, X_{b2}\right), \dots, m_{ab}\left(X_{an}, X_{bn}\right)\right)^{\intercal}.
\end{aligned}
\end{equation}

We further assume that each \(m_r\) lies in the Hilbert space \(L^2\left(\Omega\right)\), with orthonormal basis \(\{\phi_k : k \geq 0\}\), satisfying \(L^2\left(\Omega\right) = \mathrm{span}\left\{\phi_0, \phi_1, \phi_2, \dots\right\}\). 
Thus, \(m_r\) admits the expansion
\(
m_r\left(x\right) = \sum_{k=0}^{\infty} b_{rk} \phi_k\left(x\right)\) with \(b_{rk} = \int_{\Omega} m_r\left(x\right)\phi_k\left(x\right) \, \mathrm{d}x.
\)
The Hermite functions form an orthonormal basis in \(L^2(\mathbb{R})\). 
We define
\(
\phi_i(x) = \left(\sqrt{\pi} 2^i i !\right)^{-1 / 2} H_i(x) \exp\left(-\frac{x^2}{2}\right)\), 
where \(H_i(x)\) are Hermite polynomials orthogonal with respect to the density \(\exp(-x^2)\). 

Introducing a truncation parameter \(K=K_n\)\footnote{$K=K_n$ may depend on the sample size \(n\), and typically increases slowly with $n$ to balance bias and variance in the estimation.
We will discuss the choice of \(K_n\) in Section \ref{sec4}.}, we approximate the \(r^{\text{th}}\) univariate function by
\(
\bm{m}_r\left(\bm{X}_r\right) = \Phi\left(\bm{X}_r\right)\boldsymbol{b}_r + \boldsymbol{\xi}_r,
\)
where \(\boldsymbol{b}_r = \left(b_{r0}, b_{r1}, \dots, b_{r,K_n-1}\right)^{\intercal}\), and \(\Phi\left(\bm{X}_r\right)\) is the basis matrix:
\begin{equation}\label{Eqn: Basis_Matrix}
\Phi(\bm{X}_r) = \begin{bmatrix}
    \phi_0(X_{r1}) & \phi_1(X_{r1}) & \dots & \phi_{K_n-1}(X_{r1}) \\
    \phi_0(X_{r2}) & \phi_1(X_{r2}) & \dots & \phi_{K_n-1}(X_{r2}) \\
    \vdots & \vdots & \ddots & \vdots \\
    \phi_0(X_{rn}) & \phi_1(X_{rn}) & \dots & \phi_{K_n-1}(X_{rn}) \\
\end{bmatrix}.
\end{equation}
Here, \(\boldsymbol{\xi}_r\) denotes the approximation error, i.e., $\boldsymbol{\xi}_r = \sum_{k=K_n}^{\infty} b_{rk}\phi_k\left(\bm{X}_r\right)$.

For the interaction part, we similarly assume that \(m_{ab}\) lies in \(L^2\left(\Omega^2\right)\). 
It can then be expressed as
\(
m_{ab}\left(x, y\right) = \sum_{k=0}^{\infty} \sum_{d=0}^{\infty} \gamma_{abkd} \phi_k\left(x\right)\phi_d\left(y\right).
\)
Applying the truncation \(K = K_n\), we naturally obtain the following tensor product:
\begin{equation}
\footnotesize
\begin{aligned}
&\Phi\left(X_{ai}\right)\otimes \Phi\left(X_{bi}\right)
\\
&\quad =\left(\phi_0\left(X_{ai}\right)\phi_0\left(X_{bi}\right), \phi_0\left(X_{ai}\right)\phi_1\left(X_{bi}\right), \cdots, \phi_1\left(X_{ai}\right)\phi_0\left(X_{bi}\right), \cdots, \phi_{K-1}\left(X_{ai}\right)\phi_{K-1}\left(X_{bi}\right)\right).
\end{aligned}
\end{equation}

We define the basis matrix as
\begin{equation}
\Phi\left(\boldsymbol{X}_{a}\right)\otimes \Phi\left(\boldsymbol{X}_{b}\right) :=
\begin{bmatrix}
\Phi\left(X_{a1}\right)\otimes \Phi\left(X_{b1}\right) \\
\vdots \\
\Phi\left(X_{an}\right)\otimes \Phi\left(X_{bn}\right)
\end{bmatrix},
\end{equation}
which is the tensor product of the matrices \(\Phi\left(\boldsymbol{X}_{a}\right)\) and \(\Phi\left(\boldsymbol{X}_{b}\right)\) defined in \autoref{Eqn: Basis_Matrix}.
Thus, the function \(\bm{m}_{ab}\) can be approximated as
\(
\bm{m}_{ab}\left(\boldsymbol{X}_a, \boldsymbol{X}_b\right) = \Phi\left(\boldsymbol{X}_{a}\right)\otimes \Phi\left(\boldsymbol{X}_{b}\right) \cdot \boldsymbol{\gamma}_{ab} + \boldsymbol{\xi}_{ab},
\)
where: $\boldsymbol{\gamma}_{ab} = \left( \boldsymbol{\gamma}_{ab0}^{\intercal}, \boldsymbol{\gamma}_{ab1}^{\intercal}, \cdots, \boldsymbol{\gamma}_{ab,K_n-1}^{\intercal} \right)^{\intercal}$, $\boldsymbol{\gamma}_{abk} = \left(\gamma_{abk0}, \gamma_{abk1}, \cdots, \gamma_{abk,K_n-1}\right)^{\intercal}$ and \(\boldsymbol{\xi}_{ab}\) denotes the approximation error.

\rmk
Researchers may use different orthonormal bases when the Hilbert space is equipped with a density function \(\pi\), i.e.,
\(
L^2\left(\mathbb{R}; \pi\right) = \left\{ f : \int_{\mathbb{R}} f^2(x) \pi(x) \, \mathrm{d}x < \infty \right\}.
\)
For example, when \(\pi(x) = \exp\{-x^2\}\), the Hermite polynomials \(\{H_i, i \geq 0\}\), after normalization, form an orthonormal basis. Unlike the Hermite functions \(\phi_i\) that we are using, additional restrictions are necessary when taking \(H_i\) as an orthogonal basis:
(a) For any \(r\in \bm{R}\), restrict \(b_{r0} = 0\) to avoid double-counting the constant effect.
(b) For any \((a,b) \in \bm{S}\), restrict \(\gamma_{ab00} = \gamma_{ab01} = \gamma_{ab10} = 0\) to avoid counting the constant effect and the single-variable effect.

\subsection{Model Approximation and Relevant Covariates}

The true model can be approximated as
\begin{equation}
\bm{Y} = \bm{\theta} + \sum_{l \in \bm{L}} \alpha_l \bm{X}_l + \sum_{r \in \bm{R}} \bm{\Phi}(\bm{X}_r)\bm{b}_r + \sum_{(s, s') \in \bm{S}} \bm{\Phi}(\bm{X}_s) \otimes \bm{\Phi}(\bm{X}_{s'}) \bm{\gamma}_{ss'} + \bm{\Xi} + \bm{U},
\end{equation}
where \(\bm{\Xi} = \sum_{r}\bm{\xi}_r  + \sum_{a,b}\bm{\xi}_{ab}\) denotes the total approximation error, which can be small when the truncation order $K_n$ is sufficiently large, and \(\bm{U}\) is the \(n \times 1\) vector of idiosyncratic errors \(\epsilon_i\). 
We denote the non-zero coefficients by \(\bm{\alpha}_{\bm{L}}, \bm{b}_{\bm{R}}, \bm{\gamma}_{\bm{S}}\), i.e., 
\begin{equation}
\begin{aligned}
& \bm{\alpha}_{\bm{L}}=\left(\alpha_{l_1}, \cdots, \alpha_{l_L}  \right)^{\intercal} \quad \text{for } l_i \in \bm{L}, \qquad \bm{b}_{\bm{R}} = \left(\bm{b}_{r_1}^{\intercal},\cdots,\bm{b}_{r_R}^{\intercal}\right)^{\intercal} \quad \text{for } r_i \in \bm{R}, \\
& \bm{\gamma}_{\bm{S}}=\left(\bm{\gamma}_{s_1,s_1'}^{\intercal},\cdots,\bm{\gamma}_{s_S,s_S'}^{\intercal}\right)^{\intercal} \quad \text{for } (s_i, s_i') \in \bm{S}.
\end{aligned}
\end{equation}
We combine them into the single non-zero coefficient vector \(\bm{\beta}_{\bm{P}_1} = (\bm{\alpha}_{\bm{L}}^{\intercal}, \bm{b}_{\bm{R}}^{\intercal}, \bm{\gamma}_{\bm{S}}^{\intercal})^{\intercal}\). 
Let $P_1 = \mathrm{dim}(\bm{\beta}_{\bm{P}_1})$, where \(\mathrm{dim}(\cdot)\) denotes the dimension of a vector. 

To facilitate the selection procedure, we define all relevant covariates as $\bm{Z}_{P_1} = \left(\bm{Z}_{\bm{L}}, \bm{Z}_{\bm{R}}, \bm{Z}_{\bm{S}}\right)$, where
\begin{equation}\label{eqn:non0var}
\begin{aligned}
    & \bm{Z}_{\bm{L}} = \left(\bm{X}_{l_1}, \dots, \bm{X}_{l_L}\right) \quad \text{for } l_i \in \bm{L}, \\
    & \bm{Z}_{\bm{R}} = \left(\bm{\Phi}(\bm{X}_{r_1}), \dots, \bm{\Phi}(\bm{X}_{r_R})\right) \quad \text{for } r_i \in \bm{R}, \\
    & \bm{Z}_{\bm{S}} = \left(\bm{\Phi}(\bm{X}_{s_1}) \otimes \bm{\Phi}(\bm{X}_{s_1'}), \dots, \bm{\Phi}(\bm{X}_{s_S}) \otimes \bm{\Phi}(\bm{X}_{s_S'})\right) \quad \text{for } (s_i, s_i') \in \bm{S}.
\end{aligned}
\end{equation}
\section{Estimation} 
\label{sec3}

Based on the model above, our method aims to select the relevant subset \(\bm{Q}\) from \(\bm{P}\) and estimate both unknown parameters and functions. 

Our procedure builds on the seminal works of LASSO \citep{Tibshirani1996}, adaptive LASSO \citep{Zou2006}, and Group LASSO \citep{Yuan2006}, and we show that S-LASSO enjoys the \textit{oracle properties} \citep{Fan2001,Wang2008}, which require: \textbf{(1) Selection consistency}: correctly identifies relevant variables with probability approaching 1; \textbf{(2) Estimation efficiency}: yields asymptotically normal estimators for the selected variables, as if the true model were known \textit{a priori}.

\subsection{Step 1: Initial Selection}

First, we include all possible linear, nonlinear, and pairwise interaction forms of \(\bm{P}\) in the initial LASSO step, which yields
\begin{equation}\label{eq:step1}
\begin{aligned}
\min_{\bm{\alpha}, \bm{b}, \bm{\gamma}} ||\bm{Y}-\sum_{l=0}^{P}\alpha_{l}\bm{X}_{l}-&\sum_{r=1}^{P}{\Phi} (\bm{X}_{r})\bm{b}_{r}-\sum_{a=1}^{P-1}\sum_{b>a}^{P}{\Phi} (\bm{X}_{a})\otimes {\Phi}(\bm{X}_{b})\bm{\gamma} _{ab}||^{2}_{2}\\&+\lambda _n\left( \sum_{l=0}^{P}|\alpha _l|+\sum_{r=1}^{P}|\bm{b} _{r}| +\sum_{a=1}^{P-1}\sum_{b>a}^{P}|\bm{\gamma}_{ab}|\right),
\end{aligned}
\end{equation}
where the intercept is included in \(\alpha_0\), and we define \(\bm{X}_0 = \bm{1}\).
Here, \(|\cdot|\) denotes the \(l_1\) norm, \(||\cdot||_2\) denotes the \(l_2\) norm, and \(\lambda_n > 0\) is a data-driven tuning parameter. 
This step provides preliminary information for the second-stage selection.

\rmk
If prior information or conventional specifications are available, we may instead define prior candidate sets \(\bm{L}_0\), \(\bm{R}_0\), and \(\bm{S}_0\). In that case, \autoref{eq:step1} becomes
\begin{equation}\label{eq:step1turns}
\begin{aligned}
	\min_{\boldsymbol{\alpha },\boldsymbol{b},\boldsymbol{\gamma }} ||\boldsymbol{Y}-\sum_{l\in \boldsymbol{L}_0}^{}{\alpha _l}\boldsymbol{X}_l-&\sum_{r\in \boldsymbol{R}_0}^{}{\Phi}(\boldsymbol{X}_r)\boldsymbol{b}_r-\sum_{\left( a,b \right) \in \boldsymbol{S}_0}^{}{\Phi (\boldsymbol{X}_a)\otimes \Phi (\boldsymbol{X}_b)\boldsymbol{\gamma }_{ab}}||_{2}^{2}\\
	&+\lambda _n\left( \sum_{l\in \boldsymbol{L}_0}^{}{|}\alpha _l|+\sum_{r\in \boldsymbol{R}_0}^{}{|}\boldsymbol{b}_r|+\sum_{\left( a,b \right) \in \boldsymbol{S}_0}^{}{|\boldsymbol{\gamma }_{ab}|} \right) .\\
\end{aligned}
\end{equation}

\subsection{Step 2: Group Selection}

Following \cite{Zou2006} and \cite{Wang2008}, we next apply adaptive group LASSO using the estimates from Step 1. 
Specifically, we construct penalty weights and then substitute all basis functions into the adaptive group LASSO criterion:
\begin{equation}\label{eq:step2}
\begin{aligned}
\min_{\bm{\alpha}, \bm{b}, \bm{\gamma}} ||\bm{Y}-\sum_{l=0}^{P}\alpha_{l}\bm{X}_{l}-&\sum_{r=1}^{P}{\Phi} (\bm{X}_{r})\bm{b}_{r}-\sum_{a=1}^{P-1}\sum_{b>a}^{P}{\Phi} (\bm{X}_{a})\otimes {\Phi}(\bm{X}_{b})\bm{\gamma} _{ab}||^{2}_{2}\\&+\widetilde{\lambda}_n\left( \sum_{l=0}^{P}\widehat{w}_{\alpha_l}|\alpha _l|+\sum_{r=1}^{P}\widehat{w}_{\bm{b} _{r}}\left\|\bm{b} _{r}\right\|_2 +\sum_{a=1}^{P-1}\sum_{b>a}^{P}\widehat{w}_{\bm{\gamma}_{ab}}\left\|\bm{\gamma}_{ab}\right\|_2\right),
\end{aligned}
\end{equation}
where \(\widetilde{\lambda}_n\) is a tuning parameter, and the weights are defined as
\begin{equation}
\begin{aligned}
&\widehat{w}_{\alpha_l} = 
\begin{cases}
\left|\widetilde{\alpha}_l\right|^{-1}, & \text{if } \left|\widetilde{\alpha}_l\right| > 0, \\
\infty, & \text{if } \left|\widetilde{\alpha}_l\right| = 0,
\end{cases} \quad 
\widehat{w}_{\bm{b}_r} = 
\begin{cases}
\sqrt{N_R} \left\|\widetilde{\bm{b}}_r\right\|_2^{-1}, & \text{if } \left\|\widetilde{\bm{b}}_r\right\|_2 > 0, \\
\infty, & \text{if } \left\|\widetilde{\bm{b}}_r\right\|_2 = 0,
\end{cases} \\
&\widehat{w}_{\bm{\gamma}_{ab}} = 
\begin{cases}
\sqrt{N_S} \left\|\widetilde{\bm{\gamma}}_{ab}\right\|_2^{-1}, & \text{if } \left\|\widetilde{\bm{\gamma}}_{ab}\right\|_2 > 0, \\
\infty, & \text{if } \left\|\widetilde{\bm{\gamma}}_{ab}\right\|_2 = 0.
\end{cases}
\end{aligned}
\end{equation}
Here, \(\widetilde{\alpha}\), \(\widetilde{\bm{b}}\), and \(\widetilde{\bm{\gamma}}\) are the initial estimates from Step 1. 
The terms \(N_R = K_n\), and \(N_S = K_n^2\) denote the group sizes and control the strength of the penalties. The subgradient descent algorithm for S-LASSO is given in the Appendix {C}.

Apart from the relevant variables defined in \autoref{eqn:non0var}, we define the irrelevant variables as $\bm{Z}_{P_2} = \left( \bm{Z}_{\bm{L}^C}, \bm{Z}_{\bm{R}^C}, \bm{Z}_{\bm{S}^C} \right)$, with
\begin{equation}
\begin{aligned}
    & \bm{Z}_{\bm{L}^C} = \left( \bm{X}_{l_1}, \dots, \bm{X}_{l_{L^C}} \right) \quad \text{for } l_i \in \bm{L}^C, \\
    & \bm{Z}_{\bm{R}^C} = \left( \bm{\Phi}(\bm{X}_{r_1}), \dots, \bm{\Phi}(\bm{X}_{r_{R^C}}) \right) \quad \text{for } r_i \in \bm{R}^C, \\
    & \bm{Z}_{\bm{S}^C} = \left( \bm{\Phi}(\bm{X}_{s_1}) \otimes \bm{\Phi}(\bm{X}_{s_1'}), \dots, \bm{\Phi}(\bm{X}_{s_{S^C}}) \otimes \bm{\Phi}(\bm{X}_{s_{S^C}'}) \right) \quad \text{for } (s_i, s_i') \in \bm{S}^C,
\end{aligned}
\end{equation}
where $L^C:=|\bm{L}^C|$, $R^C:=|\bm{R}^C|$, $S^C:=|\bm{S}^C|$.
The \(n \times P_1\) matrix \(\bm{Z}_{P_1}\) collects all relevant components, with \(\bm{\beta}_{P_1}\) as its coefficient vector. 
Similarly, the \(n \times P_2\) matrix \(\bm{Z}_{P_2}\) collects all irrelevant components, with \(\bm{\beta}_{P_2}\) as its coefficient vector. 
Combining them, we write $\bm{Z} = \left(\bm{Z}_{P_1}, \bm{Z}_{P_2}\right)$ and $\bm{\beta}_{P_Z} = \left(\bm{\beta}_{P_1}^{\intercal},\bm{\beta}_{P_2}^{\intercal}\right)^{\intercal}$. 
The model can then be expressed as $E\left(\bm{Y} \mid \bm{P}\right) = \bm{Z} \bm{\beta}_{P_Z} + \bm{\Xi}$, where \(\bm{\Xi}\) is the approximation error, which depends on the truncation parameter \(K_n\).

The resulting model specification is
\begin{equation}
\bm{h}(\widehat{\bm{Q}}) \approx \sum_{l \in \widehat{\bm{L}}} \beta_l \bm{X}_l + \sum_{r \in \widehat{\bm{R}}} m_r(\bm{X}_r) + \sum_{(s, s') \in \widehat{\bm{S}}} m_{ss'}(\bm{X}_s, \bm{X}_{s'}),
\end{equation}
where $\widehat{\bm{Q}}$ is the set of selected variables, and $\widehat{\bm{L}} = \left\{l : |\widehat{\alpha}_l| > 0\right\}$, $\widehat{\bm{R}} = \left\{r : \left\|\widehat{\bm{b}}_r\right\|_2 > 0\right\}$, $\widehat{\bm{S}} = \left\{(s, s') : \left\|\widehat{\bm{\gamma}}_{ss'}\right\|_2 > 0\right\}$, and \(\widehat{\alpha}\), \(\widehat{\bm{b}}\), and \(\widehat{\bm{\gamma}}\) are the estimates from Step 2.
When prior candidate sets \(\bm{L}_0\), \(\bm{R}_0\), and \(\bm{S}_0\) are available, selection is restricted to these subsets: $\widehat{\bm{L}} = \left\{l \in \bm{L}_0 : |\widehat{\alpha}_l| > 0\right\}$, $\widehat{\bm{R}} = \left\{r \in \bm{R}_0 : \left\|\widehat{\bm{b}}_r\right\|_2 > 0\right\}$ and $\widehat{\bm{S}} = \left\{(s, s') \in \bm{S}_0 : \left\|\widehat{\bm{\gamma}}_{ss'}\right\|_2 > 0\right\}$.
\section{Theoretical Results} 
\label{sec4}

In this section, we establish the oracle properties of S-LASSO. 
We begin with the assumptions used in the analysis.

\subsection{Assumptions}
\label{as1}

\begin{assu}\label{as1.1}
    The samples \(\left(\boldsymbol{X}_{i}^{\intercal}, \epsilon_i\right)\) are independent and identically distributed, with \(E(\epsilon_i|\boldsymbol{X}) = 0\) and \(\mathrm{Var}(\epsilon_i|\boldsymbol{X}) = \sigma^2\). 
    Additionally, we assume \(E\left(\|\boldsymbol{X}_i\|^2\right) < \infty\). 
    Furthermore, the tail of the distribution of \(\epsilon_i\) is bounded by an exponential-type tail, i.e., there exist constants \(z_1, z_2 > 0\) such that, $
    P\left(|\epsilon_i| > x\right) \leq \exp\left\{-\left(x/z_1\right)^{z_2}\right\}.$
\end{assu}

\rmk 
These conditions are standard moment and tail assumptions in high-dimensional models. 
They are used to control the stochastic score terms that arise in the selection and estimation analysis.

For any \(K\), define the sample covariance matrix of all components by
\begin{equation}
\boldsymbol{V}^{(n)} = \frac{1}{n} (\boldsymbol{Z}_{P_1}, \boldsymbol{Z}_{P_2})^\intercal (\boldsymbol{Z}_{P_1}, \boldsymbol{Z}_{P_2}) =
\begin{bmatrix}
\boldsymbol{V}_{11}^{(n)} & \boldsymbol{V}_{12}^{(n)} \\
\boldsymbol{V}_{21}^{(n)} & \boldsymbol{V}_{22}^{(n)}
\end{bmatrix}.
\end{equation}
Since the Hermite functions are uniformly bounded: $\sup_{x \in \mathbb{R}} \sup_{k \geq 0} |\phi_k(x)| \leq \pi^{-1/4}$, \(\phi_k(X_{ri})\) and \(\phi_{k_1}(X_{ai}) \phi_{k_2}(X_{bi})\) are bounded and possess moments of all orders, so \(\bm{V}^{(n)}\) is well behaved. 
Let
\(
\boldsymbol{Z} = \left(\boldsymbol{Z}_{P_1}, \boldsymbol{Z}_{P_2}\right) = \left(Z_1, \dots, Z_{P_1}, Z_{P_1+1}, \dots, Z_{P_1+P_2}\right) = \left(\boldsymbol{z}_1^\intercal, \dots, \boldsymbol{z}_n^\intercal\right)^\intercal.
\)
Lemma {A.2} in the Appendix shows that $\frac{1}{n} Z_j^\intercal Z_j \leq M_1$ for $j = 1, 2, \dots, P_1+P_2$, for some constant \(M_1\). 
This parallels Condition (5) in \cite{zhao2006model}.

\rmk 
Assumption \autoref{as1.1} can be modified for other settings, for example, by replacing independence with suitable mixing conditions in time series applications, provided that it still implies $\frac{1}{n} Z_j^\intercal Z_j \leq M_1$.

\begin{assu}\label{as1.2}
There exist constants \(c_1\), \(c_2\), \(c_3\), \(c_4\) with \(0 \leq c_1 < c_2 < 1\) and \(c_3, c_4 > 0\), such that: 
(a)\(K_n \leq O\left(n^{\frac{c_1}{2}}\right)\);
(b) For the \(P_1 \times 1\) vector \(\bm{\beta}_{\bm{P}_1}\), let \(d\left(\bm{\beta}_{\bm{P}_1}\right)\) be the smallest element in \(|\bm{\beta}_{\bm{P}_1}|\), we suppose $    n^{\frac{1-c_2}{2}} d\left(\bm{\beta}_{\bm{P}_1}\right) \geq c_4$;
(c) \(\rho_{\min}(\boldsymbol{V}_{11}) > c_3\) uniformly for any \(n\), where \(\rho_{\min}(\cdot)\) denotes the smallest eigenvalue.
\end{assu}

\rmk 
Condition (a) controls the truncation level, condition (b) imposes a minimum signal strength, and condition (c) excludes severe multicollinearity by requiring the smallest eigenvalue of \(\bm{V}_{11}^{(n)}\) to be bounded away from zero.

To clarify condition (a), \cite{zhao2006model} assume \(P_1 = O\left(n^{c_1}\right)\) and \(P_2 = o\left(n^{(c_2-c_1)k}\right)\) for some integer \(k > 0\), provided \(E\left(\epsilon_i^{2k}\right) < \infty\). 
In our setting, \(P_1 = O\left(K_n^2\right)\) and \(P_2 = O\left(P^2 \cdot K_n^2\right)\). 
Under condition (a), these growth rates are compatible with \cite{zhao2006model}; see Lemma {A.4} in the Appendix.

\begin{assu}\label{as1.4}
The sample covariance matrix \(\bm{V}^{(n)}\) satisfies the strong irrepresentable condition of \citet{zhao2006model}: there exists a positive \(P_{2}\times 1\) vector \(\bm{\eta}\) such that
\begin{equation}\label{eqn:irrecon}
\left|\bm{V}^{(n)}_{21}(\bm{V}^{(n)}_{11})^{-1}\mathrm{sign}(\bm{\beta}_{\bm{P}_1})\right|
\leq \bm{1}-\bm{\eta},
\end{equation}
where the inequality holds element-wise.
\end{assu}

\rmk
Assumption~\ref{as1.4} is needed for the sign consistency in Theorem \ref{thm1.1}.
The weak irrepresentable condition,
\(
\left|\bm{V}^{(n)}_{21}(\bm{V}^{(n)}_{11})^{-1}\mathrm{sign}(\bm{\beta}_{\bm{P}_1})\right|
\leq \bm{1},
\)
is less restrictive but insufficient for the sign consistency result.

\begin{assu}\label{as1.5}
The functions \(m_r(x)\) and \(m_{ab}(x, y)\) satisfy the following smoothness conditions:
(a) For any \(r \in \bm{R}\), \(m_r(x)\) is differentiable up to order \(l_1\) (\(l_1 \geq 1\)), and for \(i = 0, 1, \dots, l_1\), \(m_r^{(l_1-i)}(x) x^i \in L^2(\Omega)\). 
We denote this class as \(F^{l_1}(\Omega)\) for simplicity;
(b) For any \((a, b) \in \bm{S}\), the functions \(f_{abx_0}(y) = m_{ab}(x_0, y)\) and \(f_{aby_0}(x) = m_{ab}(x, y_0)\) belong to \(F^{l_2}(\Omega)\) for all \(x_0, y_0\), with \(l_2 \geq 3\).
\end{assu}

\rmk 
These smoothness conditions are in line with \cite{DGT2016AoS} and ensure the convergence of the orthogonal series approximation.

\subsection{Selection Consistency}

The following definition formalizes the sign consistency.
For brevity, we write \(\widehat{\bm{\beta}} =_s \bm{\beta}\) for \(\mathrm{sign}(\widehat{\bm{\beta}}) = \mathrm{sign}(\bm{\beta})\), where $\mathrm{sign}(\cdot)$ maps positive entries to 1, negative entries to -1 and zeros to 0.

\begin{defn}
Let \(\widehat{\bm{\beta}}\) be the estimator of \(\bm{\beta}\) through the S-LASSO. Then, \(\widehat{\bm{\beta}}\) is \textit{Sign Consistent} with \(\bm{\beta}\) if and only if there exists a sequence of tuning parameters \(\lambda_n\) (or \(\widetilde{\lambda}_n\)) such that
\(
\lim \limits_{n \to \infty} P\left(\widehat{\bm{\beta}}\left(\lambda_n\right) =_s \bm{\beta}\right) = 1\).
\end{defn}

Recall $P_1 = \dim(\bm{\beta}_{\bm{P}_1})$, $P_2  = \dim(\bm{\beta}_{\bm{P}_2})$, and let $P_Z = P_1 + P_2$ denote the dimension of $\bm{\beta}_{\bm{P}_Z}$. The following theorem establishes selection consistency.

\begin{thm}\label{thm1.1}[Selection Consistency]
Under Assumptions \autoref{as1.1}-\autoref{as1.5},
for any \(\lambda_n\) satisfying \(\frac{\lambda_n}{\sqrt{n}}= o\left(n^{\frac{c_2-c_1}{2}}\right)\) and \(\frac{\lambda_n}{n^{\frac{1}{2}+\delta}} \to \infty\) for some \(\delta > 0\), if \(K_n \to \infty\) as \(n \to \infty\), then the estimator from the Step 1 are sign consistent. 
Specifically, $\lim \limits_{n \to \infty}  P\left( \widetilde{\bm \beta}_{\bm P_Z} \neq_s \bm \beta_{\bm P_Z} \right) = 0$.
\end{thm}

\rmk
Theorem \ref{thm1.1} shows that Step 1 correctly identifies the relevant components and shrinks the irrelevant ones to zero. Because the adaptive weights in \autoref{eq:step2} assign an infinite penalty to zero first-step estimates, the same zero pattern is preserved in Step 2. Hence,
\(
\lim \limits_{n \to \infty}  P\left( \widehat{\bm \beta}_{\bm P_Z} \neq_s \bm \beta_{\bm P_Z} \right) = 0.
\)

\rmk 
In classical LASSO theory \citep{zhao2006model}, the irrepresentable condition is almost necessary and sufficient for exact sign consistency. 
However, it can be restrictive because all entries in \autoref{eqn:irrecon} are required to hold simultaneously. 
This is particularly demanding in asset pricing applications, where firm characteristics are often correlated \citep{freyberger2017dissecting}. 
When the condition fails, the global sign consistency event $\{\widehat{\bm{\beta}}_{\bm{P}_Z} = _s \bm{\beta}_{\bm{P}_Z}\}$ may fail as well.
However, if only a small fraction of components violate the irrepresentable inequality, then the entrywise sign consistency probability for a typical component can still remain high.
Appendix {D.3} reports an additional diagnostic analysis under empirically calibrated dependence structures.

\subsection{Asymptotic Normality}

We now study the consistency and asymptotic normality of the estimator in Step 2.

\begin{thm}\label{thm1.2}[Consistency]
Under Assumptions \autoref{as1.1}-\autoref{as1.5}, if \(\widetilde{\lambda}_n = o(n^{c_2/2})\) and \(K_n = O(n^\nu)\) with \(0 < \nu < 1/4\), then:
(a) $\widehat{\bm{\alpha}}_{\bm{L}}$ is consistent, i.e., \(\widehat{\bm{\alpha}}_{\bm{L}} \xrightarrow{P} \bm{\alpha}_{\bm{L}}\).
(b) $\widehat{m}_r$ is $L_2$-consistent, i.e., $\left\|m_r - \widehat{m}_r  \right\|_{L^2} = o_P\left(1\right)$, for any $r \in \bm R$.
(c) $\widehat{m}_{ab}$ is $L_2$-consistent, i.e., $\left\|m_{ab} - \widehat{m}_{ab}  \right\|_{L^2} = o_P\left(1\right)$, for any $\left(a,b\right) \in \bm S$.
(d) Specifically, for the rate of convergence, we have
$
\left\| m_r-\widehat{m}_r \right\| _{L^2}\le O_P\left( n^{-\min \left\{ \frac{1}{2}-\frac{\nu}{2},\frac{l_1\nu}{2} \right\}} \right)$, $\left\| m_{ab}-\widehat{m}_{ab} \right\| _{L^2}^{}\le O_P\left( n^{-\min \left\{ \frac{1}{2}-\nu ,\frac{(l_2-1)\nu}{2} \right\}} \right) .
$
\end{thm}

\rmk
This theorem shows that, once irrelevant variables are removed, Step 2 yields consistent estimators. The $L_2$-consistency is standard in nonparametric statistics; see \cite{Huang2010}.

We next study the asymptotic distribution of the estimators. Define
\begin{equation}
\bm{\mathcal{T}}\left(\bm{x}\right) = \left( \begin{array}{c}
	\widehat{\boldsymbol{\alpha }}_{\boldsymbol{L}}-\boldsymbol{\alpha }_{\boldsymbol{L}},\\
	\left( \widehat{m}_{r_j}\left( x_{r_j} \right) -m_{r_j}\left( x_{r_j} \right) \right) _{1\le j\le R}\\
	\left( \widehat{m}_{a_jb_j}\left( x_{a_j},x_{b_j} \right) -m_{a_jb_j}\left( x_{a_j},x_{b_j} \right) \right) _{1\le j\le S}\\
\end{array} \right),
\end{equation}
which collects all parameters and functions at a given $\bm{x} \in \Omega^Q$. Its dimension is fixed because $L+R+S\le 2Q+\frac{Q(Q-1)}{2}$. The next theorem shows that \(\bm{\mathcal{T}}(\bm{x})\) is asymptotically normal.

\begin{thm}\label{thm1.3}[\(\sqrt{n}\)-Consistency]
Under the conditions in Theorem \autoref{thm1.2} and assume
$
K_n = O(n^{\nu})$ with $\max\left\{\frac{6}{6l_1 - 5}, \frac{6}{6l_2 - 17}\right\} \leq \nu < 1/4,
$
then we have
\begin{equation}
\sqrt{n}\left( \sigma ^2\boldsymbol{G}\left( \boldsymbol{x} \right) \boldsymbol{V}_{11}^{-1}\boldsymbol{G}\left( \boldsymbol{x} \right) ^{\intercal} \right) ^{-1/2} \bm{\mathcal{T}}(\bm{x}) \xrightarrow{d}N\left( 0,\bm I_{L+R+S} \right),
\end{equation}
where 
$
\boldsymbol{G}(\bm{x})= \mathrm{diag}\left\{\bm{I}_L, \mathrm{diag}\left\{ \left( \Phi \left( x_{r_j} \right) \right) _{1\le j\le R} \right\}, \mathrm{diag}\left\{ \left( \Phi \left( x_{a_j} \right) \otimes \Phi \left( x_{b_j} \right) \right) _{1\le j\le S} \right\} \right\}.
$
\end{thm}

\rmk
Theorem \ref{thm1.3} implies that the linear coefficients in Step 2 are \(\sqrt{n}\)-consistent and that, for any given \(\bm{x} \in \Omega^{Q}\), the estimated nonlinear and interactive functions are also \(\sqrt{n}\)-consistent. This establishes the oracle property of S-LASSO.

Based on Theorem \autoref{thm1.3}, for a specific \(r_0\), we have $\frac{\sqrt{n}\left( \widehat{m}_{r_0}\left( x \right) -m_{r_0}\left( x \right) \right)}{\sqrt{\Phi \left( x \right) \mathbf{\Sigma }_{\boldsymbol{b}_{r_0}}^{-1}\Phi \left( x \right) ^{\intercal}}}\xrightarrow{d} N\left( 0,\sigma^2 \right)$, where 
\(
\mathbf{\Sigma}_{\bm{b}_{r_0}} := \mathrm{Cov}(\Phi(X_{r_0 1}) \mid \bm{z}_{-r_0, 1}),
\)
and 
\begin{equation}
\bm{z}_{-r_0, 1} = \left(\{X_{l1} : l \in \bm{L}\}, \{\Phi(X_{r1}) : r \neq r_0, r \in \bm{R}\}, \{\Phi(X_{a1}) \otimes \Phi(X_{b1}) : (a, b) \in \bm{S}\}\right).
\end{equation}
Similarly, for the interactive function $m_{a_0,b_0}(x,y)$, we have
\begin{equation}
\frac{\sqrt{n}\left( \widehat{m}_{a_0b_0}\left( x,y \right) -m_{a_0b_0}\left( x,y \right) \right)}{\sqrt{\left( \Phi \left( x \right) \otimes \Phi \left( y \right) \right) \mathbf{\Sigma }_{\boldsymbol{\gamma }_{a_0b_0}}^{-1}\left( \Phi \left( x \right) \otimes \Phi \left( y \right) \right) ^{\intercal}}}\rightarrow N\left( 0,\sigma^2 \right),
\end{equation}
where $\sigma^2=\mathrm{Var}(\epsilon_i|\bm{X})$, and we define
\(
\mathbf{\Sigma}_{\bm{\gamma}_{a_0 b_0}} := \mathrm{Cov}(\Phi(X_{a_0,1}) \otimes \Phi(X_{b_0,1}) \mid \bm{z}_{-(a_0, b_0), 1})
\)
and
\begin{equation}
\footnotesize
\bm{z}_{-(a_0, b_0), 1} = \left(\{X_{l1} : l \in \bm{L}\}, \{\Phi(X_{r1}) : r \in \bm{R}\}, \{\Phi(X_{a1}) \otimes \Phi(X_{b1}) : (a, b) \neq (a_0, b_0), (a, b) \in \bm{S}\}\right).
\end{equation}

\section{Simulation Study}
\label{sec5}
\subsection{Data Generating Processes}

In this section, we examine the finite-sample performance of S-LASSO using Monte Carlo simulations.
We generate the dataset
\(\{(y_i,X_i^{\intercal})^{\intercal}: i=1,\cdots,n\}\) from
\begin{equation}
y_i = \beta x_{1i} + m_1(x_{2i}) + m_2(x_{3i}, x_{4i}) + \epsilon_i, \quad i = 1, \dots, n,
\end{equation}
which contains one linear effect, one nonlinear univariate effect, and one interaction effect.
The covariates satisfy \(X_i \sim N(\bm{0},\bm{\Sigma})\), where each marginal distribution is \(N(0,1)\), and the errors satisfy \(\epsilon_i \sim N(0,\sigma^2)\).
The simulation design is specified as follows.

\begin{enumerate}

    \item \textbf{Covariate correlation matrix.}
    The correlation matrix \(\bm{\Sigma}=(\sigma_{ij})_{P\times P}\) is constructed from an empirical distribution calibrated from the real data. 
    Specifically, each off-diagonal entry \(\sigma_{ij}\) is drawn from the empirical distribution of pairwise correlations among the 102 firm characteristics used in Section~\ref{sec6}. 
    The correlation structure of empirical data is discussed in Appendix {D.3}.

    \item \textbf{Two designs for the true functions.}
    We consider two data-generating processes.
    In DGP 1, we set \(\beta=1\), \(m_1(x_{2i})=2.24\cos(x_{2i})\), and \(m_2(x_{3i},x_{4i})=1.43\sin(x_{3i}+x_{4i})\).
    In DGP 2, we set \(\beta=1\), \(m_1(x_{2i})=2.24\cos(x_{2i})\), and \(m_2(x_{3i},x_{4i})=1.91\sin(x_{3i}x_{4i})\).
    The constants normalize the signal components to have approximately unit variance. 
    DGP 1 represents a genuinely bivariate interaction, whereas DGP 2 is intentionally favorable to the FNW benchmark because the interaction depends on the product \(x_{3i}x_{4i}\).

    \item \textbf{Sample size and error level.}
    We fix the sample size at \(n=200\) and consider two dimensions, \(P=100\) and \(P=300\), corresponding to \(P<n\) and \(P>n\).
    We set \(\sigma=1\), which corresponds to a unit noise-to-signal ratio.
\end{enumerate}

\subsection{Simulation Results}

We compare the proposed S-LASSO with several benchmark models:

\begin{enumerate}
    \item \textbf{LASSO.} 
    We apply the standard LASSO to the original linear covariates only. 

    \item \textbf{LN-LASSO.} 
    We apply the LASSO to a dictionary containing linear terms and univariate nonlinear basis expansions, but no interaction terms. 

    \item \textbf{FNW.} 
    We implement the approach of \citet{freyberger2017dissecting}, which includes linear terms, nonlinear univariate effects, and product-type interaction terms taking the form \(m(x_a x_b)\).
\end{enumerate}

The key benchmark is the FNW method of \citet{freyberger2017dissecting}. 
To provide a balanced comparison, DGP 2 is intentionally specified in a way that is favorable to this method. 
Since the true functions are smooth, we use Hermite basis expansions for nonlinear terms and tensor-product expansions for interaction terms, with truncation level \(K=n^{1/6}\).
In the simulation, we allow all potential pairwise interactions to enter the candidate dictionary.

For all methods and steps, the tuning parameter is selected by the Bayesian Information Criterion (BIC),
defined as \(\mathrm{BIC} = \log\!\left(\mathrm{MSE}\right) + df \cdot \frac{\log n}{n}\),
where \(\mathrm{MSE}\) denotes the mean squared error and \(df\) is the number of nonzero coefficients.
Each setting is repeated 100 times, and we report the means and standard deviations in \autoref{tab:simu_res}.

\begin{table}[h!]
\centering
\caption{
Simulation Results.
\footnotesize{
This table compares S-LASSO with three benchmarks: standard linear LASSO, LN-LASSO with linear and univariate nonlinear terms but no interactions, and FNW with univariate nonlinear terms and product-form interactions. 
`LIN.', `NONL.', and `INTE.' report correct detection rates for the true linear, nonlinear, and interaction effects, respectively. 
`FP' denotes the false positive rate and `MSE' denotes mean squared error. 
Entries are averages over 100 replications, with standard deviations in parentheses.}}
\label{tab:simu_res} 
\setlength{\tabcolsep}{6.75pt}
\setstretch{1.5}
\scriptsize
\begin{tabular}{llccccccccccc}
\toprule
\multirow{2}{*}{\textbf{Model}} &  & \multicolumn{4}{c}{\textbf{Selection ($P=100$)}} & \multirow{2}{*}{\textbf{MSE}} & \multicolumn{1}{l}{} & \multicolumn{4}{c}{\textbf{Selection ($P=300$)}} & \multirow{2}{*}{\textbf{MSE}} \\
\cmidrule(lr){3-6} \cmidrule(lr){9-12}
                                &  & LIN.       & NONL.     & INTE.     & FP        &                               & \multicolumn{1}{l}{} & LIN.     & NONL.     & INTE.      & FP         &                               \\ \midrule
\rowcolor{gray!25}\multicolumn{13}{l}{\textbf{$\bullet$ Panel A: DGP 1}}                                                                                                                                                                       \\
\multirow{2}{*}{LASSO}          &  & 1          & 0         & 0         & 0.04      & 2.6090                        &                      & 1        & 0         & 0          & 0.01       & 2.6633                        \\
                                &  & (0)        & (0)       & (0)       & (0.02)    & (0.3006)                      &                      & (0)      & (0)       & (0)        & (0.01)     & (0.3429)                      \\
\multirow{2}{*}{LN-LASSO}       &  & 0.99       & 1         & 0         & 0.36      & 0.0081                        &                      & 1        & 1         & 0          & 0.16       & 0.0009                        \\
                                &  & (0.10)     & (0)       & (0)       & (0.01)    & (0.0013)                      &                      & (0)      & (0)       & (0)        & (0.00)     & (0.0001)                      \\
\multirow{2}{*}{FNW}            &  & 0.99       & 1         & 0.53      & 0.12      & 0.0068                        &                      & 1        & 1         & 0.15       & 0.05       & 0.0134                        \\
                                &  & (0.10)     & (0)       & (0.50)    & (0.00)    & (0.0018)                      &                      & (0)      & (0)       & (0.36)     & (0.00)     & (0.0040)                      \\
\multirow{2}{*}{S-LASSO}        &  & 1          & 1         & 1         & 0.12      & 0.0062                        &                      & 1        & 1         & 0.87       & 0.05       & 0.0126                        \\
                                &  & (0)        & (0)       & (0)       & (0.00)    & (0.0015)                      &                      & (0)      & (0)       & (0.34)     & (0.00)     & (0.0036)                      \\[1ex]
\rowcolor{gray!25}\multicolumn{13}{l}{\textbf{$\bullet$ Panel B: DGP 2}}                                                                                                                                                                       \\
\multirow{2}{*}{LASSO}          &  & 1          & 0         & 0         & 0.01      & 3.1231                        &                      & 1        & 0         & 0          & 0.00       & 3.0995                        \\
                                &  & (0)        & (0)       & (0)       & (0.01)    & (0.3382)                      &                      & (0)      & (0)       & (0)        & (0.00)     & (0.3223)                      \\
\multirow{2}{*}{LN-LASSO}       &  & 0.98       & 1         & 0         & 0.36      & 0.0295                        &                      & 1        & 1         & 0          & 0.16       & 0.0009                        \\
                                &  & (0.14)     & (0)       & (0)       & (0.04)    & (0.2105)                      &                      & (0)      & (0)       & (0)        & (0.00)     & (0.0001)                      \\
\multirow{2}{*}{FNW}            &  & 1          & 1         & 1         & 0.12      & 0.0069                        &                      & 1        & 1         & 1          & 0.05       & 0.0175                        \\
                                &  & (0)        & (0)       & (0)       & (0.00)    & (0.0017)                      &                      & (0)      & (0)       & (0)        & (0.00)     & (0.0048)                      \\
\multirow{2}{*}{S-LASSO}        &  & 1          & 1         & 1         & 0.12      & 0.0075                        &                      & 1        & 1         & 1          & 0.05       & 0.0203                        \\
                                &  & (0)        & (0)       & (0)       & (0.00)    & (0.0017)                      &                      & (0)      & (0)       & (0)        & (0.00)     & (0.0065)                      \\ \bottomrule
\end{tabular}
\end{table}

\autoref{tab:simu_res} shows four main patterns. 
First, S-LASSO recovers the linear and nonlinear effects almost perfectly across both DGPs and both values of \(P\), and allowing all pairwise interactions does not weaken the recovery of the univariate effect. 
Second, under DGP~2, which favors FNW, both FNW and S-LASSO achieve perfect selection, with FNW having a slightly smaller MSE. 
Third, under DGP~1, where the interaction is not of product form, S-LASSO has a clear advantage: when \(P=300\), FNW identifies the interaction in only \(15\%\) of replications, whereas S-LASSO maintains an interaction recovery rate of \(87\%\). 
Finally, increasing \(P\) mainly affects interaction recovery, as expected from the growing number of candidate interactions, while linear and nonlinear recovery remain stable.

\begin{figure}[h!]
\centering
  \begin{subfigure}[b]{0.75\textwidth}
    \includegraphics[width=1\textwidth]{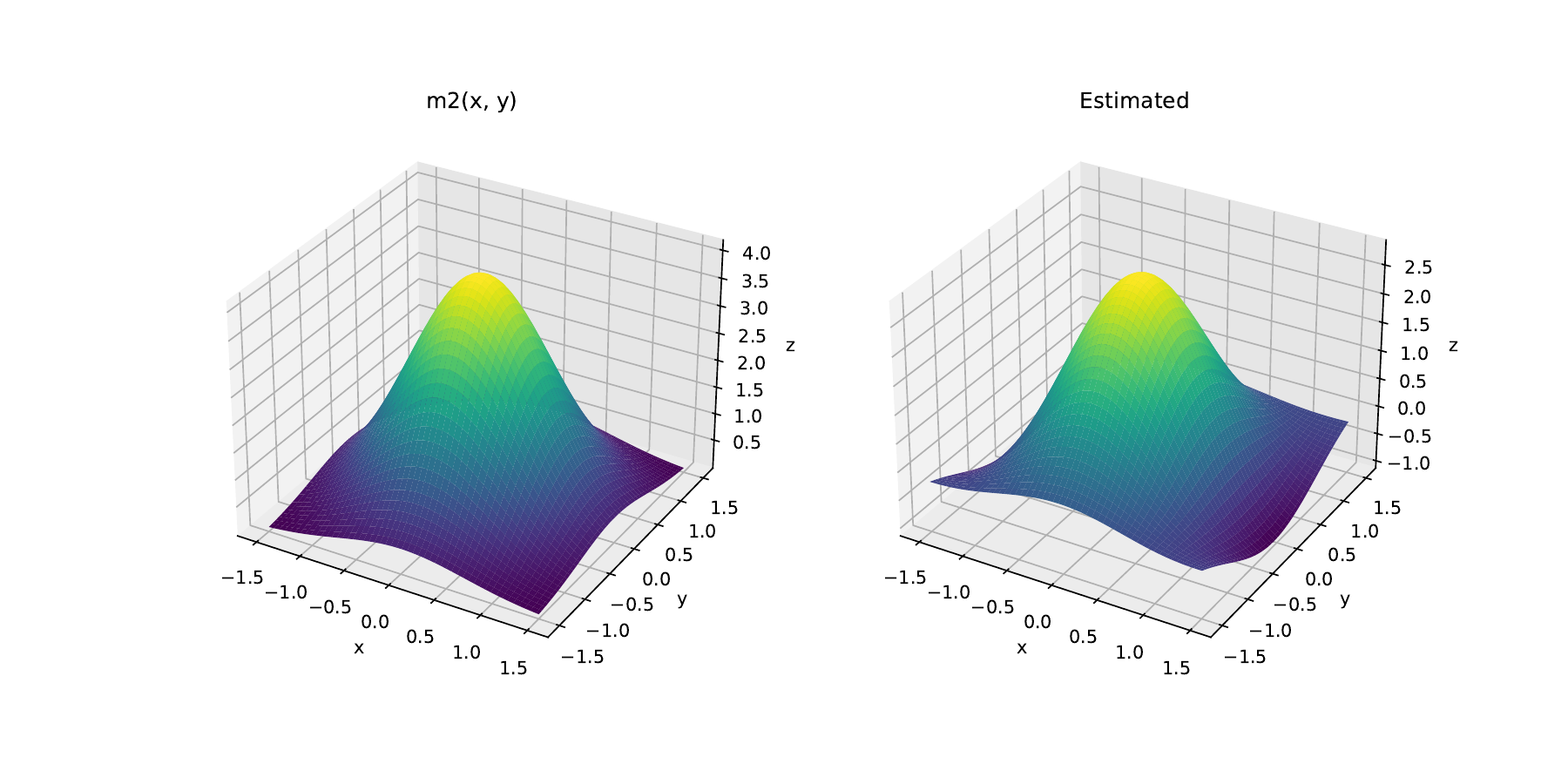}
    \caption{Interaction Function}
  \end{subfigure}
  
  \begin{subfigure}[b]{1\textwidth}
    \centering
    \includegraphics[width=0.52\textwidth]{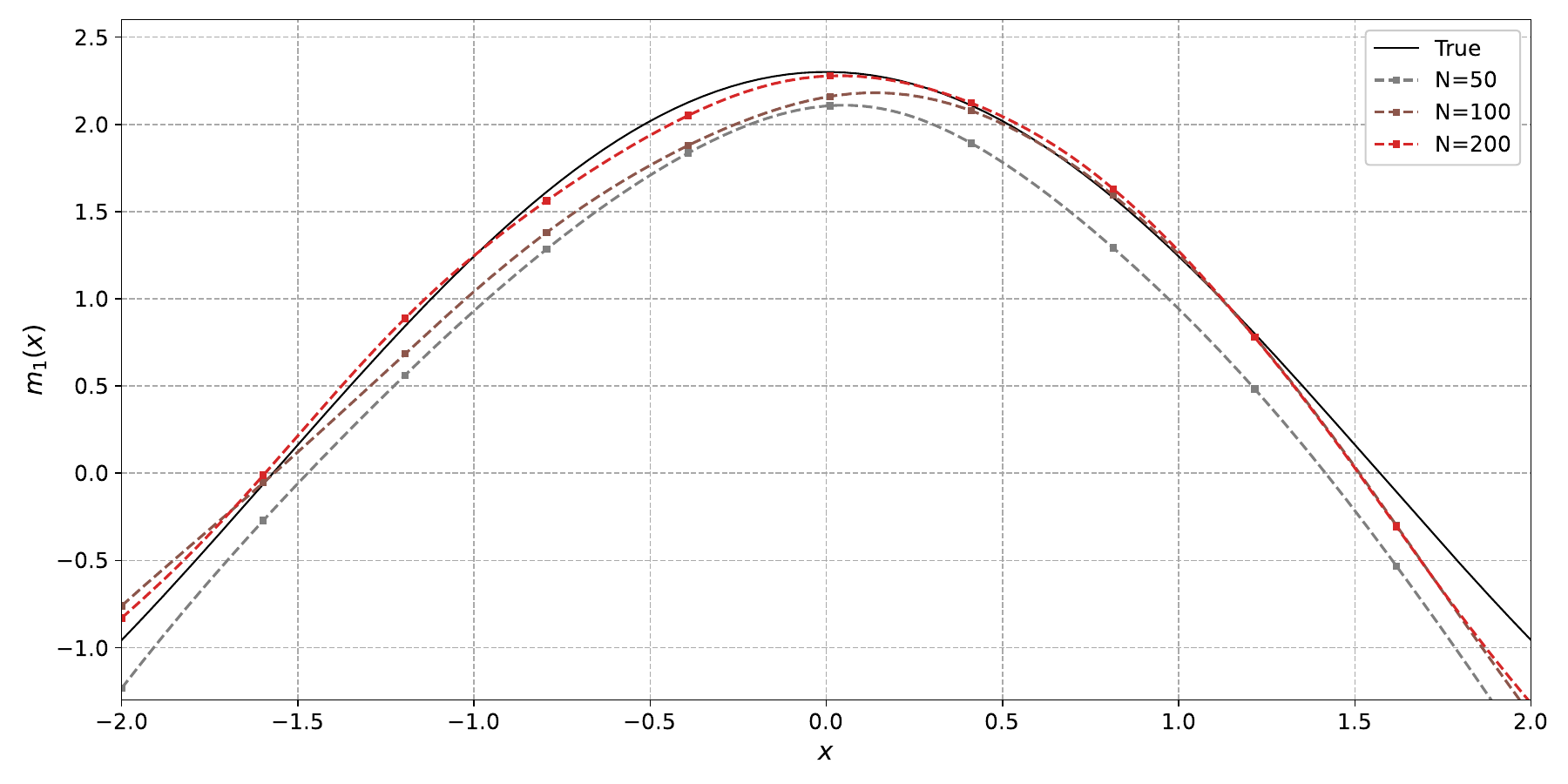}
    \caption{Nonlinear Function}
  \end{subfigure}  
  \caption{Estimated and True Functions.}
  \label{fig:funcs}
\end{figure}

\autoref{fig:funcs} reports the estimated functions under DGP~1 with \(n=200\) and \(P=300\). 
Both the nonlinear and interaction estimates closely track the true functions.

\section{Empirical Study}
\label{sec6}

In this section, we revisit the question in \citet{freyberger2017dissecting} of how firm characteristics predict stock returns in a nonparametric framework. They model asset returns as additive nonparametric functions of lagged firm characteristics and allow pairwise multiplicative interactions with firm size:
\begin{equation}\label{8}
E(\bm{Y}|\bm{P}) = \theta \boldsymbol{1} + \sum_{r=1}^{R} m_{r}(\bm{X}_{r}) + \sum_{s \neq \text{size}} m_{s,\mathrm{size}}(\bm{X}_{s} \cdot \bm{X}_{\text{size}}),
\end{equation}
where the unknown interaction functions take inputs of the form 
$
\bm{X}_{s} \cdot \bm{X}_{\text{size}}=(X_{s1}X_{\mathrm{size},1},\cdots,X_{sn}X_{\mathrm{size},n})^{\intercal}.
$
This specification restricts interactions to pairwise product-form interactions and excludes linear terms.

We apply S-LASSO to re-examine the effects of firm characteristics on stock returns and their interactions with firm size. Our specification is
\begin{equation}\label{9}
E(\bm{Y}|\bm{P}) = \theta \boldsymbol{1} + \sum_{l=1}^{L} \alpha_l \boldsymbol{X}_l + \sum_{r=1}^{R} m_{r}(\bm{X}_{r}) + \sum_{s \neq \text{size}} m_{s,\mathrm{size}}(\bm{X}_{s}, \bm{X}_{\text{size}}),
\end{equation}
where the notation is consistent with \autoref{8}. Relative to \citet{freyberger2017dissecting}, our model adds linear terms and allows more flexible interactions with firm size. We also assume that both slope coefficients and characteristic functions are time-invariant within each rolling window.

\rmk
The interaction terms in S-LASSO can be chosen according to the empirical question. 
In the main application, we focus on interactions with firm size, motivated by the central role of size in asset pricing \citep{HongLimStein2000JF} and by the benchmark of \citet{freyberger2017dissecting}. 
These size-centered interactions improve interpretability and facilitate comparison with the existing nonparametric benchmark. 
Appendix {D.2} further reports a full-specification check using all pairwise interactions over consecutive four-year subsamples; firm size is selected in all ten subsamples and is the most frequently connected characteristic, supporting the main empirical specification.

\subsection{Data}

Following \citet{freyberger2017dissecting}, our empirical objective is to identify which firm characteristics predict monthly returns and through which functional forms. We use the 102 characteristics listed in \citet{GHZ2017RFS}, all of which can be constructed from CRSP, Compustat, and/or I/B/E/S. 

\begin{table}[h!]
\centering
\caption{
Summary of Data.
\footnotesize{This table summarizes the return data and characteristics data used in this paper.
For the return data, the descriptive statistics are calculated across all sample dates and sample firms, where p(J-B) denotes the $p$-value for Jarque-Bera test.}}
\label{tab: data}
\tabcolsep=0.18cm
\footnotesize
\setstretch{1.5}
\begin{tabular}{cccccccccccc}
\toprule
\textbf{Variable} & \textbf{Mean} & \textbf{Std.} & \textbf{t-stat.} & \textbf{Skew.} & \textbf{Kurt.} & \textbf{p (J-B)} & \textbf{Q1} & \textbf{Median} & \textbf{Q3} & \textbf{Count} \\ \midrule
Return            & 1.18\%        & 19.10\%       & 0.06             & 8.28           & 443.82         & 0.00             & -6.67\%     & 0.00\%          & 7.29\%      & 2198914  \\ \bottomrule 
\end{tabular}
\vspace{0.5cm}

\scriptsize
\tabcolsep=0.115cm
\setstretch{1.25}
\begin{tabular}{>{\itshape}p{2.5cm}p{1cm}p{1cm}p{0.1cm}>{\itshape}p{2.5cm}p{1cm}p{1cm}p{0.1cm}>{\itshape}p{2.5cm}p{1.2cm}p{1cm}p{0.1cm}}
\toprule
\rm{\textbf{Variable}} & \textbf{Mean} & \textbf{Std.} & \textbf{} & \rm{\textbf{Variable}} & \textbf{Mean} & \textbf{Std.} & \textbf{} & \rm{\textbf{Variable}} & \textbf{Mean} & \textbf{Std.} \\ \midrule
absacc            & 0.090         & 0.123         &           & grltnoa           & 0.087         & 0.177         &           & rd\_sale          & 0.445         & 5.023         \\
acc               & -0.041        & 0.195         &           & herf              & 0.079         & 0.088         &           & realestate        & 0.244         & 0.238         \\
aeavol            & 0.855         & 2.254         &           & hire              & 0.090         & 0.372         &           & retvol            & 0.034         & 0.027         \\
age               & 11.911        & 9.766         &           & idiovol           & 0.067         & 0.039         &           & roaq              & -0.005        & 0.059         \\
agr               & 0.162         & 0.473         &           & ill               & 0.000         & 0.000         &           & roavol            & 0.029         & 0.049         \\
baspread          & 0.055         & 0.067         &           & invest            & 0.072         & 0.214         &           & roeq              & -0.004        & 0.162         \\
beta              & 1.098         & 0.679         &           & ipo               & 0.067         & 0.251         &           & roic              & -0.141        & 1.121         \\
betasq            & 1.659         & 1.917         &           & lev               & 2.219         & 4.717         &           & rsup              & 0.005         & 0.248         \\
bm                & 0.689         & 0.640         &           & lgr               & 0.259         & 0.845         &           & salecash          & 54.950        & 169.899       \\
cash              & 0.160         & 0.218         &           & maxret            & 0.077         & 0.077         &           & saleinv           & 26.351        & 74.817        \\
cashdebt          & -0.042        & 1.326         &           & mom12m            & 0.134         & 0.652         &           & salerec           & 10.910        & 59.060        \\
cashpr            & -1.298        & 58.048        &           & mom1m             & 0.010         & 0.158         &           & secured           & 0.496         & 0.540         \\
cfp               & 0.043         & 0.732         &           & mom36m            & 0.341         & 1.146         &           & securedind        & 0.470         & 0.499         \\
chcsho            & 0.112         & 0.328         &           & mom6m             & 0.054         & 0.383         &           & sfe               & -0.367        & 6.623         \\
chfeps            & -0.001        & 0.428         &           & ms                & 3.650         & 2.011         &           & sgr               & 0.191         & 0.626         \\
chinv             & 0.008         & 0.090         &           & mve               & 12.002        & 2.299         &           & sin               & 0.009         & 0.093         \\
chmom             & -0.004        & 0.597         &           & nanalyst          & 4.153         & 6.476         &           & sp                & 2.029         & 3.401         \\
chnanalyst        & -0.020        & 1.415         &           & nincr             & 0.948         & 1.412         &           & std\_dolvol       & 0.851         & 0.428         \\
chtx              & 0.001         & 0.013         &           & operprof          & 0.773         & 1.229         &           & std\_turn         & 4.561         & 11.947        \\
cinvest           & -0.009        & 2.828         &           & orgcap            & 0.009         & 0.014         &           & stdacc            & 3.863         & 32.052        \\
convind           & 0.121         & 0.327         &           & pchcurrat         & 0.081         & 0.604         &           & stdcf             & 9.721         & 87.271        \\
currat            & 3.294         & 5.104         &           & pchdepr           & 0.114         & 0.541         &           & sue               & -0.005        & 0.149         \\
depr              & 0.285         & 0.433         &           & pchgm\_pchsale    & -0.079        & 1.015         &           & tang              & 0.532         & 0.173         \\
disp              & 0.158         & 0.506         &           & pchquick          & 0.113         & 0.727         &           & tb                & -0.133        & 1.976         \\
divi              & 0.029         & 0.168         &           & pchsale\_pchinvt  & -0.050        & 0.838         &           & turn              & 1.177         & 1.894         \\
divo              & 0.030         & 0.173         &           & pchsale\_pchrect  & -0.049        & 0.721         &           & zerotrade         & 1.238         & 3.326         \\
dolvol            & 11.435        & 3.115         &           & pchsale\_pchxsga  & 0.025         & 0.453         &           & bm\_ia            & 25.902        & 758.801       \\
dy                & 0.015         & 0.033         &           & pchsaleinv        & 0.176         & 1.051         &           & cfp\_ia           & 7.120         & 434.171       \\
ear               & 0.003         & 0.084         &           & pctacc            & -1.025        & 7.770         &           & chatoia           & -0.010        & 0.268         \\
egr               & 0.149         & 0.753         &           & pricedelay        & 0.150         & 1.031         &           & chempia           & -0.136        & 0.747         \\
ep                & -0.041        & 0.366         &           & ps                & 4.238         & 1.817         &           & chpmia            & 0.333         & 14.107        \\
fgr5yr            & 12.830        & 15.511        &           & quick             & 2.600         & 4.478         &           & pchcapx\_ia       & 8.669         & 85.143        \\
gma               & 0.351         & 0.370         &           & rd                & 0.134         & 0.370         &           & indmom            & 0.144         & 0.302         \\
grcapx            & 0.959         & 3.472         &           & rd\_mve           & 0.057         & 0.131         &           & mve\_ia           & -124.156      & 8388.043    \\ \bottomrule
\end{tabular}
\end{table}

We also collect monthly stock returns from CRSP and adjust delisting returns following \citet{hou2015digesting}. The sample spans January 1980 to December 2020. 
Summary statistics are reported in \autoref{tab: data}.
Variable definitions are reported in Appendix {D.1}.
Since only four characteristics have complete coverage \citep{GHZ2017RFS}, we handle missing data using the method of \citet{FHNW2024RFS}, which applies conditional mean imputation and weighted least squares within a GMM framework.

\subsection{Specification Results}

\subsubsection{Selected Characteristics}

We apply S-LASSO to \autoref{9} to determine which characteristics predict returns and through which functional forms, with tuning parameters selected by BIC. 
To accommodate potential time-varying parameters and functional forms, we divide the full sample into ten subperiods: 1980--1983, 1984--1987, $\cdots$, 2016--2019.
The results are summarized in \autoref{tab: spe_res_short}.

\begin{table}[h!]
\centering
\caption{
Summary of Specified Characteristics by S-LASSO.
\footnotesize{
This table reports how frequently each firm characteristic is selected by S-LASSO across ten four-year subperiods. 
Panels A and B list characteristics selected frequently as linear and nonlinear effects, respectively. 
Panel C lists characteristics selected as having interaction effects with firm size, measured by market equity (`mve'). 
The tuning parameter is selected by BIC.}}
\label{tab: spe_res_short}
\scriptsize
\tabcolsep=0.45cm
\setstretch{1.45}
\begin{tabular}{lll}
\toprule
\textbf{\textit{\#Selected}} & \textbf{} & \textbf{Characteristics}                                                             \\ \midrule
\rowcolor{gray!25}\multicolumn{3}{l}{\textbf{$\bullet$ Panel A: Linear Effects}}                                                                                   \\
\textit{\#Selected} = 10     &           & cash, mom1m, turn                                                                    \\
\textit{\#Selected} = 9      &           & baspread, indmom, secured                                                            \\
\textit{\#Selected} = 8      &           & chpmia, idiovol, std\_dolvol, std\_turn, chempia, ipo, rd\_mve                       \\
\textit{\#Selected} = 7      &           & mom12m, mom36m, retvol, sue, mve\_ia, age, agr, mom6m, realestate, betasq, zerotrade \\[1ex]
\rowcolor{gray!25}\multicolumn{3}{l}{\textbf{$\bullet$ Panel B: Nonlinear Effects}}                                                                                          \\
\textit{\#Selected} = 10     &           & mom1m, indmom, mve                                                                   \\
\textit{\#Selected} = 9      &           & mom12m, chmom                                                                        \\
\textit{\#Selected} = 8      &           & std\_turn, mom36m, nincr                                                             \\
\textit{\#Selected} = 7      &           & sue, mom6m, maxret                                                                   \\[1ex]
\rowcolor{gray!25}\multicolumn{3}{l}{\textbf{$\bullet$ Panel C: Interactive Effects with `mve'}}                                                                                      \\
\textit{\#Selected} = 4      &           & mom36m, mve\_ia                                                                      \\
\textit{\#Selected} = 3      &           &                                                                                      \\
\textit{\#Selected} = 2      &           & beta                                                                                 \\
\textit{\#Selected} = 1      &           & mom1m, std\_turn, agr, nanalyst, ps, gma, ill, divi, ms, securedind      \\
\bottomrule
\end{tabular}
\end{table}

Several characteristics, such as `cash' and `turn', are selected repeatedly as linear effects, while others, especially `indmom' and `mom1m', are repeatedly selected as nonlinear effects. Overall, 99 characteristics are selected at least once in the linear component and 87 at least once in the nonlinear component, indicating substantial functional-form heterogeneity. Panel C reports interactions with the market value of equity (`mve'). In total, 13 characteristics exhibit interaction effects with size, including `mom36m', `beta', `mom1m', and `ms'.

These results yield three main findings:
\textbf{(i) Functional-Form Heterogeneity}: Firm characteristics predict stock returns through different functional forms. 
Both linear and nonlinear effects are pervasive, suggesting that restricting all characteristics to enter linearly would miss important nonlinear predictive patterns.
\textbf{(ii) Size-Dependent Effects}: We find strong interactions between firm size (`mve') and 13 characteristics, including momentum (`mom36m' and `mom1m'), volatility-related measures (`std\_turn'), and profitability indicators (`gma').
\textbf{(iii) Comparison with Existing Studies}: Relative to \cite{freyberger2017dissecting}, our specification reveals both linear effects and richer size-dependent interactions. In particular, variables such as `mom1m', `cash', and `turn' enter linearly in our models, while characteristics including `beta', `gma', and `ms' show flexible interactions with size beyond the multiplicative form.

\subsubsection{Economic Interpretations}

To interpret in-sample selection results, we visualize representative estimated functions from the 2016--2019 window.

\textbf{Representative Nonlinear Effects.}

\begin{figure}[h!]
\centering
\includegraphics[width=0.75\textwidth]{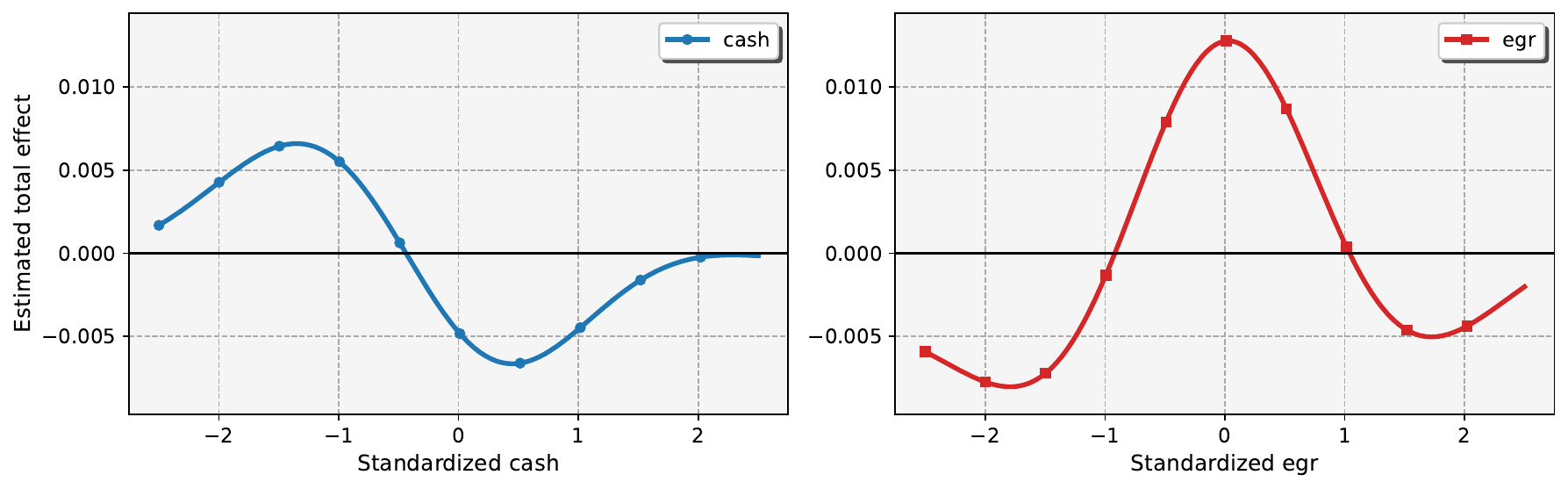}
\caption{Estimated nonlinear effects of `cash' and `egr'. 
Each curve shows the estimated total effect after centering. The horizontal axes are standardized.}
\label{fig:cash-egr-effects}
\end{figure}

We first examine two selected nonlinear effects: ‘cash’ and ‘egr’. \autoref{fig:cash-egr-effects} reports their estimated total effects after centering. For ‘cash’, the estimated function displays a non-monotonic, approximately Z-shaped pattern: firms with very low cash holdings are associated with higher predicted returns, the effect becomes negative for firms with intermediate cash holdings, and it rises again toward the upper tail of the cash distribution. This pattern suggests that cash holdings do not enter the return prediction model through a single monotone channel. One interpretation is that very low cash holdings may proxy for liquidity pressure or financing constraints, whereas very high cash holdings may reflect precautionary savings by firms with greater exposure to aggregate risk, consistent with the mechanism in \cite{Palazzo2012JFE}. For ‘egr’, growth in common shareholder equity, the estimated function is hump-shaped: the effect is negative in both tails of the distribution and positive around the middle range. This pattern indicates that common-equity growth does not enter the return prediction model through a simple linear channel. The negative effect for high equity growth is consistent with prior evidence that equity issuance, external financing, and firm expansion are informative for future returns \citep{PontiffWoodgate2008JF,CooperGulenSchill2008JF}.

\textbf{Momentum Effects and Size Interaction.}

\begin{figure}[h!]
\centering
\includegraphics[width=0.75\textwidth]{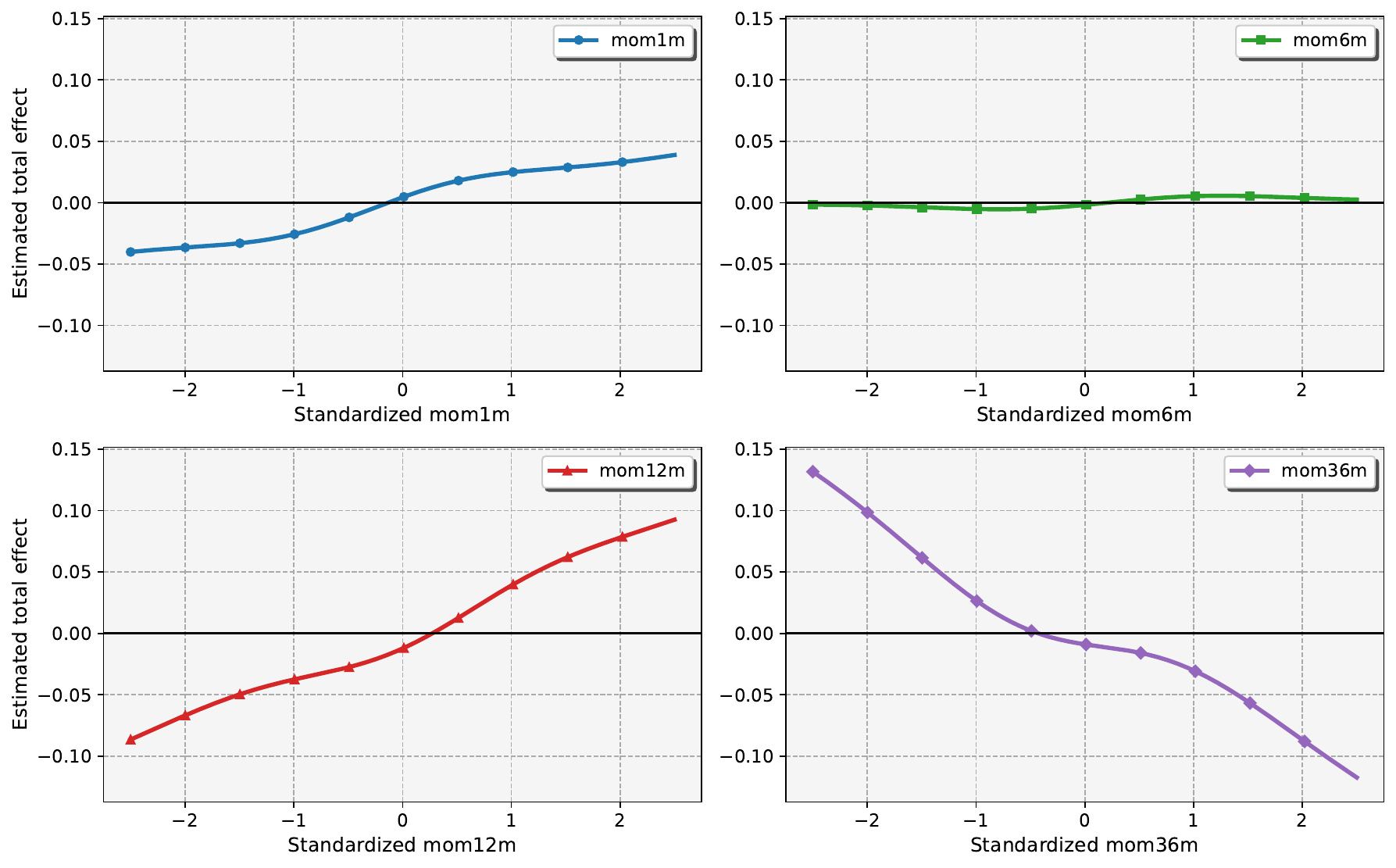}
\caption{Estimated total effects of momentum variables at different horizons. 
Each curve combines the selected linear and nonlinear components. The horizontal axis is standardized.}
\label{fig:momentum-total-effects}
\end{figure}

In the same subsample, S-LASSO assigns both linear and nonlinear components to `mom1m', `mom6m', `mom12m', and `mom36m', and selects an interaction between `mom1m' and `mve'. 
\autoref{fig:momentum-total-effects} shows that the selected momentum effects are not linear. 
The one-month momentum effect is steeper for low recent returns and becomes flatter as recent returns increase. Additionally, the six-month momentum 'mom6m' is relatively weak.
The effect of `mom36m' is negative over most of its range, while its slope also changes with the level of past returns. 
Thus, momentum is not a homogeneous predictor; its contribution depends on the horizon over which past returns are measured, consistent with evidence that short-, intermediate-, and long-horizon past returns capture different return-predictability effects \citep{NOVYMARX2012JFE}.

\begin{figure}[h!]
\centering
\includegraphics[width=0.55\textwidth]{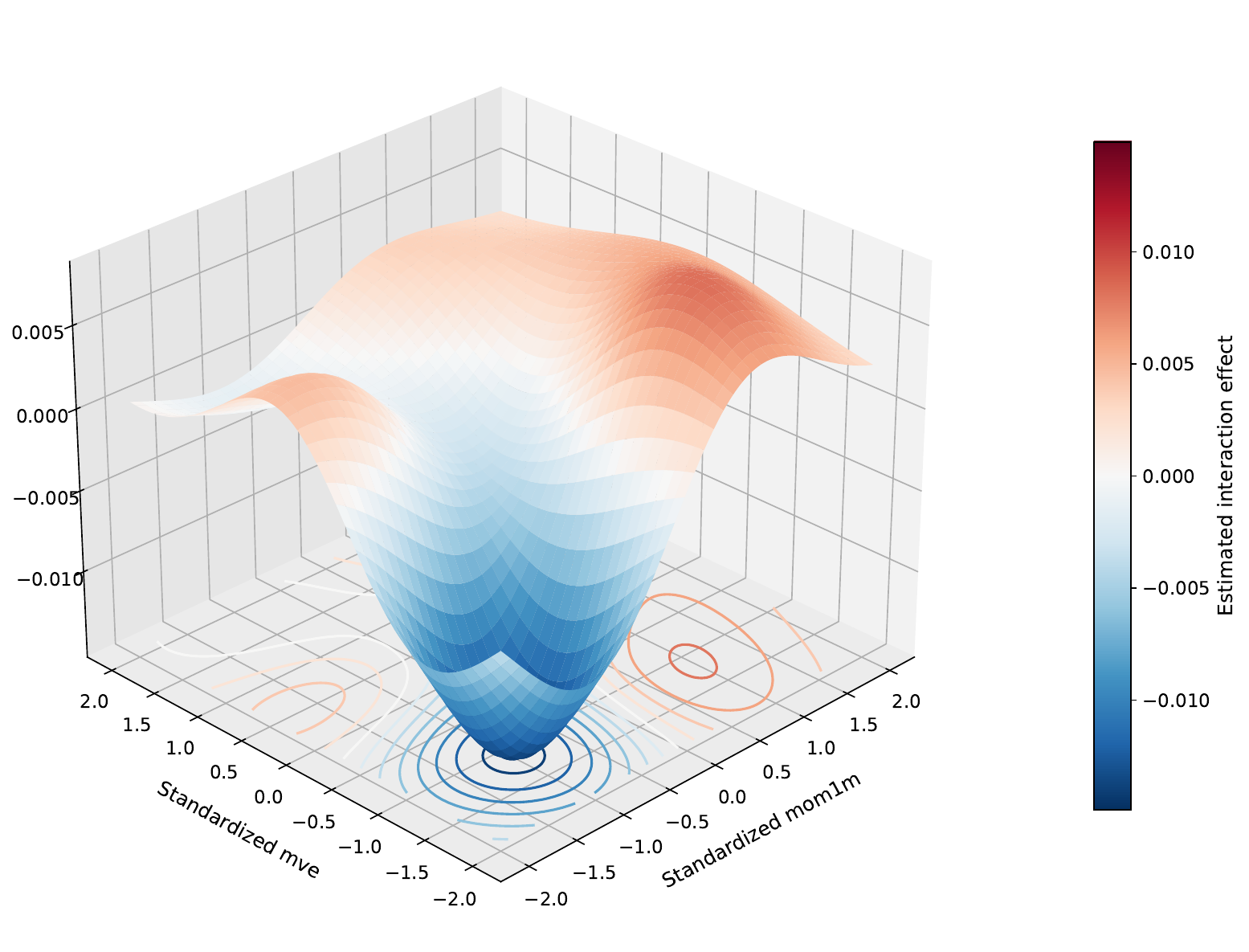}
\caption{Estimated interaction surface between one-month momentum `mom1m' and firm size `mve'. 
Both axes are standardized.}
\label{fig:mom1m-mve-interaction}
\end{figure}

\autoref{fig:mom1m-mve-interaction} plots the estimated interaction surface between one-month momentum and firm size.
The surface is non-flat: the interaction effect changes more sharply in the low-`mve' region and becomes flatter as `mve' increases. The size-dependent interaction pattern is also broadly consistent with the asset-pricing literature showing that the profitability of momentum-related strategies depends on firm size and related information frictions \citep{HongLimStein2000JF}. More broadly, classic studies suggest that return continuation and reversal patterns can vary systematically with firm size. \cite{jegadeesh1990evidence} documents the existence of short-horizon return reversals and examines serial-correlation patterns across different size groups of stocks. For longer-horizon reversal (`mom36m'), there is also evidence that the effect is substantially stronger for smaller firms than for larger firms (\cite{chopra1992measuring}).

More recent studies likewise show that many anomalies are more pronounced among small firms, including \cite{hou2017replicating} and \cite{GHZ2017RFS}. Taken together, this evidence motivates our choice to interact characteristics with firm size.

\subsection{Out-of-sample Performance}
We compare S-LASSO with benchmark methods using mean squared prediction error and long-short portfolio performance.

\begin{enumerate}
    \item \textbf{Out-of-sample Mean Squared Error}: 
    Using a rolling four-year estimation window, we estimate the model, forecast one-month-ahead returns, and compute the mean squared prediction error (MSE). Repeating this procedure through 2020 allows us to assess predictive stability across different market conditions.
    
    \item \textbf{Long--short Portfolios}: 
    Following \citet{GHZ2017RFS}, each month we rank stocks by model-predicted returns using only the preceding four years of data, and form an equal-weighted portfolio that goes long the top 20\% and short the bottom 20\% of stocks.
\end{enumerate}

\begin{table}[h!]
\centering
\caption{
Out-of-sample Performance.
\footnotesize{
The two portfolio blocks report performance over different evaluation periods, while the MSE block summarizes monthly prediction errors.}}
\label{tab:oosstats}
\scriptsize
\tabcolsep=0.125cm
\setstretch{1.45}
\begin{tabular}{ccccclcccclccc}
\toprule
\multirow{2}{*}{\textbf{Method}} & \multicolumn{4}{c}{\textbf{Portfolios (from 1984)}} &                      & \multicolumn{4}{c}{\textbf{Portfolios (from 1996)}} &  & \multicolumn{3}{c}{\textbf{MSE}}                    \\ \cline{2-5} \cline{7-10} \cline{12-14} 
                                 & Mean       & Std.       & Sharpe   & Drawdown &                      & Mean       & Std.       & Sharpe   & Drawdown   &  & Mean                 & Std.                 & Max.                 \\ \midrule
\textbf{S-LASSO}                  & 22.70\%    & 13.44\%    & 1.69     & 40.08\% & \multicolumn{1}{c}{} & 18.06\%    & 14.72\%    & 1.23     & 40.08\%    &  & 4.03                & 16.67               & 102.54              \\
\textbf{FNW}                     & 20.82\%    & 13.34\%    & 1.56     & 40.29\% & \multicolumn{1}{c}{} & 15.48\%    & 14.49\%    & 1.07     & 40.29\%    &  & 4.84                & 21.72               & 133.30              \\
\textbf{Equal Weight}            & 15.04\%    & 13.14\%    & 1.14     & 53.66\% &                      & 15.27\%    & 13.71\%    & 1.11     & 53.66\%    &  & \multicolumn{1}{l}{} & \multicolumn{1}{l}{} & \multicolumn{1}{l}{} \\ \bottomrule
\end{tabular}
\end{table}

\begin{figure}[h!]
    \centering
    \includegraphics[width=0.72\textwidth]{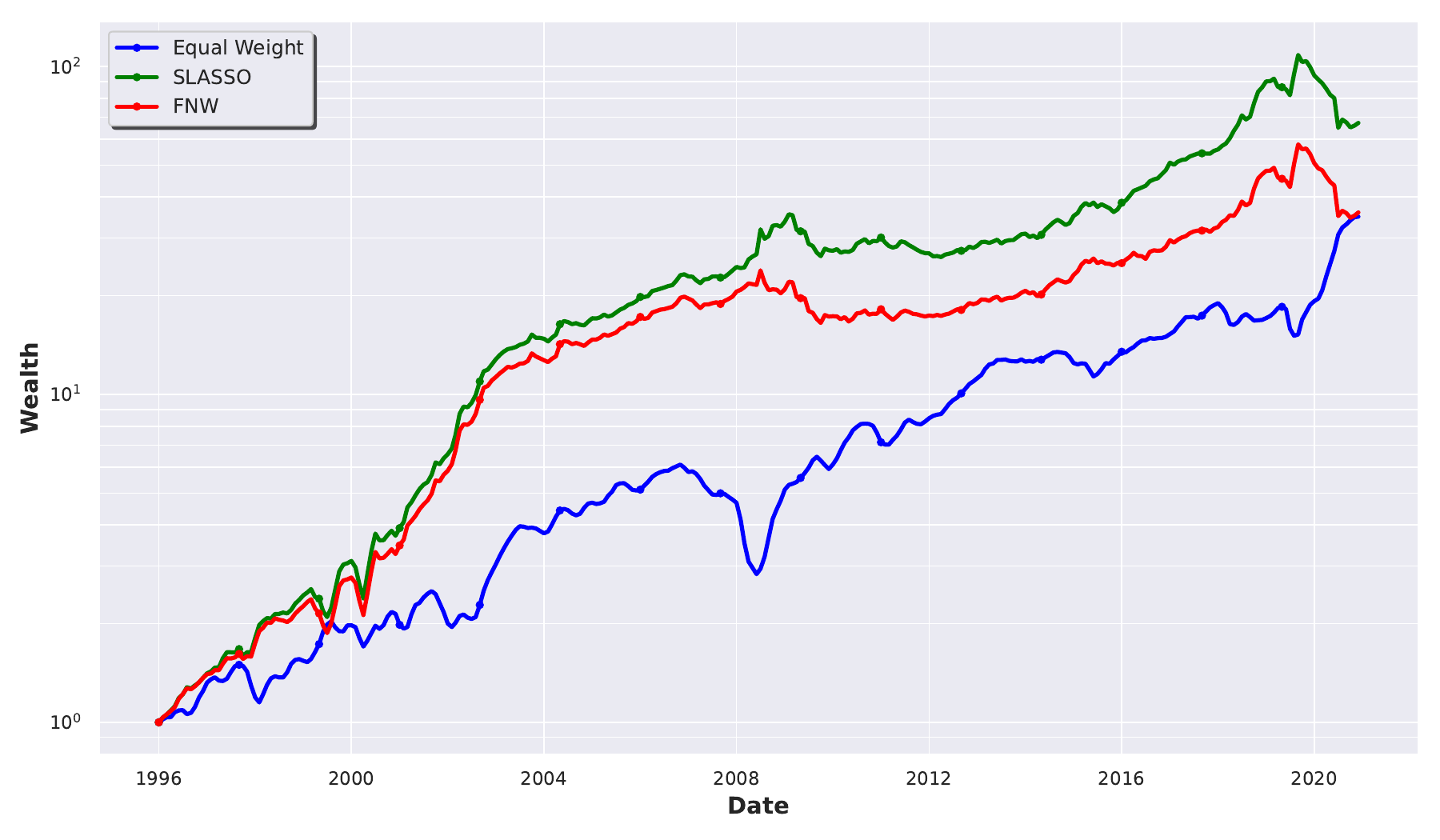}
    \caption{Out-of-sample portfolios.}
    \label{fig:oos-portfolios}
\end{figure}

\autoref{tab:oosstats} reports the out-of-sample results. 
S-LASSO outperforms FNW across the reported prediction-error and portfolio-performance metrics. 
Its mean squared prediction error falls to 4.03, a 16.7\% improvement over FNW's 4.84, and its forecasts are more stable, with an MSE standard deviation of 16.67 versus 21.72. 
The corresponding long--short portfolio earns 22.7\% per year with a Sharpe ratio of 1.69, compared with 20.8\% and 1.56 for FNW, while the maximum drawdown remains nearly identical (40.1\% versus 40.3\%). 
\autoref{fig:oos-portfolios} plots the log cumulative portfolio returns since 1996 and shows that S-LASSO consistently outperforms the benchmark strategies over the sample period.

\section{Conclusion}
\label{sec7}

This paper proposes Specification LASSO (S-LASSO), a unified framework for variable selection and model specification in high-dimensional semiparametric models. 
By combining linear, nonlinear, and pairwise interaction components, S-LASSO selects both relevant covariates and their functional forms. Theoretically, we establish oracle properties, including selection consistency and asymptotic normality. 
Monte Carlo simulations show that S-LASSO remains robust under correlated covariates and can distinguish linear, nonlinear, and interaction effects in high-dimensional settings. Empirically, we apply S-LASSO to characteristics-based asset pricing. 
The results show that firm characteristics affect returns through heterogeneous channels, including linear effects, nonlinear effects, and size-dependent interactions. 
Compared with the benchmark model of \cite{freyberger2017dissecting}, S-LASSO improves the out-of-sample predictive performance. Overall, S-LASSO offers an interpretable and practical tool for empirical researchers working with rich covariate spaces. 
Future work may extend the framework to higher-order interactions and time-varying effects.

\section*{Acknowledgments}
The authors thank the Editor, the Associate Editor, and the anonymous referee for their constructive comments and suggestions. 
Chaohua Dong acknowledges the financial support from National Natural Science Foundation of China (project code No. 72473156); Shuyi Ge is supported financially by National Natural Science Foundation of China (project code No. 72673092, 72303119, 72495155); Shaoran Li acknowledges the financial support from National Natural Science Foundation of China (project code No. 72573004, 72203007) and School of Economics, Peking University; Wen Su gratefully acknowledges support from the EPSRC Centre for Doctoral Training in Mathematics of Random Systems (EP/S023925/1). The authors report there are no competing interests to declare.


\clearpage
\pagenumbering{arabic}
\renewcommand{\thepage}{A-\arabic{page}}
\etocdepthtag{appendix}
\pdfbookmark[0]{Supplementary Appendices}{supplementary-appendices}
\begin{center}
  {\Large\bfseries Supplementary Appendices for\\[0.5em]
  ``Specification LASSO and a Flexible Characteristics-Based Asset Pricing Model ''}
\end{center}
\begingroup
\spacingset{1}
\makeatletter
\let\savedipfilectr\the@ipfilectr
\makeatother
\tableofcontents
\makeatletter
\global\let\the@ipfilectr\savedipfilectr
\makeatother
\endgroup

\appendix
\SLASSOAppendixtrue
\setcounter{equation}{0}
\setcounter{table}{0}
\setcounter{figure}{0}
\setcounter{footnote}{0}
\renewcommand*{\theHequation}{appendix.\arabic{equation}}
\renewcommand*{\theHtable}{appendix.\arabic{table}}
\renewcommand*{\theHfigure}{appendix.\arabic{figure}}
\renewcommand*{\theHfootnote}{appendix.\arabic{footnote}}
\def\ge{\geqslant}
\def\geq{\geqslant}
\def\le{\leqslant}
\def\leq{\leqslant}
\renewcommand{\theequation}{A-\arabic{equation}}
\label{appendix}

\section{Technical Lemmas}

This section establishes the fundamental properties of the basis expansions and the resulting design matrix that are used throughout the theoretical analysis. 
The first lemma records the approximation accuracy of the Hermite expansion, while the second lemma controls the column norms of the expanded design matrix. 
These two ingredients are used repeatedly below: the approximation bound controls the sieve truncation error, and the column-norm bound ensures that the stochastic score terms in the high-dimensional LASSO analysis are well behaved.

\begin{lem}\label{lemA1}
Suppose $g(x)$ is differentiable on $\mathbb{R}$ up to order $m$ ($m \geq 1$), and $g^{(m-l)}(x)x^l \in L^2(\mathbb{R})$ for $l = 0, 1, \cdots, m$. 
For the expansion $g(x) = \sum_{k=0}^{\infty} b_k \phi_k(x)$, the following results hold:\\
(1) $\sup_x \left| g\left( x \right) -\sum_{k=0}^{K-1}{b_k\phi _k\left( x \right)} \right|=o\left( 1 \right) K^{-\frac{\left( m-1 \right)}{2}-\frac{1}{12}}$;\\
(2) $\left\| g\left( x \right) -\sum_{k=0}^{K-1}{b_k\phi _k\left( x \right)} \right\| _{L^2}=o\left( 1 \right) K^{-\frac{m}{2}}$.
\end{lem}

\prf See Lemma C.1 in the Supplemental Material of \cite{DGT2016AoS}.
\qed

\begin{lem}\label{lemZ}
Let $\boldsymbol{Z} = \left( z_{ij} \right)_{n \times (P_1 + P_2)} = \left( Z_1, Z_2, \dots, Z_{P_1}, Z_{P_1+1}, \dots, Z_{P_1+P_2} \right)$.
If 
\begin{equation}
\frac{1}{n} X_k^{\intercal} X_k \leq M_0, \quad \text{for any } k,
\end{equation} 
then 
\begin{equation}
\frac{1}{n} Z_j^{\intercal} Z_j \leq M_1, \quad j = 1, 2, \dots, P_1 + P_2,
\end{equation}
where \(M_1\) is a constant.
\end{lem}

\begin{rmk}
The condition $\frac{1}{n} X_k^{\intercal} X_k \leq M_0$ is not restrictive, as it can always be satisfied by normalizing the covariates. 
The role of this condition is to put the original linear covariates and the basis-expanded covariates on a comparable scale. 
After the expansion, every column of $\boldsymbol{Z}$ is either an original covariate, a univariate basis function, or a tensor-product basis function. 
Thus, it is enough to control these three cases separately.
\end{rmk}

\prf 
We consider three cases for $z_{ij}$, corresponding to the three types of columns that can appear in the expanded design matrix $\boldsymbol{Z}$.

(a) If $z_{ij}$ is a linear term, then $Z_j$ coincides with one of the original covariate columns. 
Therefore, by the assumed normalization,
\begin{equation}
   \frac{1}{n}Z_{j}^{\intercal}Z_j = \frac{1}{n}\sum_{i=1}^n z_{ij}^{2} \leq M_0.
\end{equation}

(b) If $z_{ij}$ is generated by a nonlinear expansion, then $z_{ij}=\phi_k(X_{ri})$ for some covariate $X_r$ and some basis index $k$. 
Using the uniform boundedness of the Hermite basis functions, $|\phi_k(x)| \leq \pi^{-1/4}$, we have $z_{ij}^2 \leq \pi^{-1/2}$ uniformly in $i$. 
Hence,
\begin{equation}
   \frac{1}{n}Z_{j}^{\intercal}Z_j = \frac{1}{n}\sum_{i=1}^n z_{ij}^{2} \leq \pi^{-1/2}.
\end{equation}

(c) If $z_{ij}$ is generated by an interactive expansion, then $z_{ij}=\phi_{k_1}(X_{ai})\phi_{k_2}(X_{bi})$ for some pair $(a,b)$ and some basis indices $k_1$ and $k_2$. 
The same uniform boundedness gives $|z_{ij}| \leq \pi^{-1/2}$, and therefore $z_{ij}^2 \leq \pi^{-1}$ uniformly in $i$. 
Thus,
\begin{equation}
   \frac{1}{n}Z_{j}^{\intercal}Z_j = \frac{1}{n}\sum_{i=1}^n z_{ij}^{2} \leq \pi^{-1}.
\end{equation}

Let $M_1 = \max\{\pi^{-1/2}, M_0\}$. 
Since $\pi^{-1} \leq \pi^{-1/2}$, the three cases above are all covered by this choice of $M_1$. 
Combining the preceding bounds yields
\begin{equation}
\frac{1}{n}Z_{j}^{\intercal}Z_j \leq M_1, \quad \text{for all } j,
\end{equation}
which proves the claim.
\qed

\begin{lem}\label{lemA2}
Define $\bm W^{(n)} = \frac{1}{\sqrt{n}} \bm Z^{\intercal} \left(\bm U + \bm \Xi\right)$, and assume that Assumptions 1-4 hold. Provided $K_{n} \to \infty$, the following results hold:\\
(a) $\boldsymbol{W}^{(n)}=\boldsymbol{W}\cdot \{ 1+ o_P(1)\}$, where $\bm W\sim N\left(0,\sigma^2 \bm V\right)$, here $\boldsymbol{W}^{(n)}=\boldsymbol{W}\cdot \{ 1+ o_P(1)\}$ is defined to represent: for any $\bm{\vartheta} \in \mathbb{R}^{P_Z}$ with finite $\bm{\vartheta}^{\intercal} \bm{V} \bm{\vartheta}$, one has $\bm{\vartheta}^{\intercal}\bm{W}^{(n)} = \bm{\vartheta}^{\intercal}\bm{W} + o_P(1)$;\\
(b) Let $\boldsymbol{W}^{(n)}_1$ be the first $P_1$ entries of $\boldsymbol{W}^{(n)}$, we have $\boldsymbol{W}^{(n)}_1 =\boldsymbol{W}_1\cdot \{ 1+ o_P(1)\}$, which similarly means for any $\bm{\vartheta} \in \mathbb{R}^{P_1}$ with finite $\bm{\vartheta}^{\intercal} \bm{V}_{11} \bm{\vartheta}$, one has $\bm{\vartheta}^{\intercal}\bm{W}_1^{(n)} = \bm{\vartheta}^{\intercal}\bm{W}_1 + o_P(1)$.
\end{lem}

\prf
(a) We first separate the stochastic error part from the sieve approximation error. 
Define 
\begin{equation}
\widetilde{\bm W}^{(n)} 
= \frac{1}{\sqrt{n}}\bm Z^{\intercal}\bm U
= \frac{1}{\sqrt{n}}\sum_{i=1}^n \bm z_i^{\intercal}\varepsilon_i .
\end{equation}
For any \(K_n\), the Hermite basis functions are uniformly bounded, and hence the entries of \(\bm Z\) generated by the basis expansions have moments of all orders. 
Together with Assumption 1, this allows us to apply the central limit theorem to any finite-variance linear projection of the score vector. 
Thus, in the projected sense stated in the lemma,
\begin{equation}
\widetilde{\bm W}^{(n)}
=
\bm W\cdot \{1+o_P(1)\},
\end{equation}
where \(\bm W\sim N(0,\sigma^2\bm V)\).

It remains to show that replacing \(\bm U\) by \(\bm U+\bm\Xi\) does not change the leading stochastic order under the same projected interpretation. 
For \(\bm W^{(n)}\), define \(w_j^{(n)}\) and \(\widetilde w_j^{(n)}\) to be the \(j\)-th component of \(\bm W^{(n)}\) and \(\widetilde{\bm W}^{(n)}\), respectively. 
Then
\begin{equation}
\widetilde w_j^{(n)}-w_j^{(n)}
=
\frac{1}{\sqrt n}\sum_{i=1}^n z_{ij}\varepsilon_i
-
\frac{1}{\sqrt n}\sum_{i=1}^n z_{ij}\left(\varepsilon_i+\Xi_i\right)
=
-\frac{1}{\sqrt n}\sum_{i=1}^n z_{ij}\Xi_i,
\end{equation}
where \(\Xi_i\) is the approximation error corresponding to the \(i\)-th sample. 
The sign in the last display is immaterial for the variance bound below, but the display makes explicit that the only additional term is the sieve approximation error.

We now control this approximation term. 
Using the triangle inequality and the fact that the number of truly relevant components is fixed, we obtain
\begin{equation}\label{eqn:Xileq}
\Xi_i^2
\le
C_Q\left(
\sum_{r\in\bm R}\xi_{ir}^2
+
\sum_{(a,b)\in\bm S}\xi_{iab}^2
\right),
\quad i=1,2,\cdots,n,
\end{equation}
where
\begin{equation}
\xi_{ir}
=
m_r(X_{ir})-\Phi(X_{ir})\bm b_r
=
\sum_{k\ge K} b_{rk}\phi_k(X_{ir}),
\end{equation}
and
\begin{equation}
\xi_{iab}
=
m_{ab}(X_{ia},X_{ib})
-
\Phi(X_{ia})\otimes\Phi(X_{ib})\bm\gamma_{ab}.
\end{equation}

By Lemma \autoref{lemA1}, the univariate approximation error satisfies 
\(\xi_{ir}=o(1)K^{-\frac{l_1-1}{2}-\frac{1}{12}}\), and hence 
\(\xi_{ir}^2=o(1)K^{-l_1+\frac{5}{6}}\). 
For the bivariate approximation error, we apply the same Hermite approximation argument successively in the two coordinates. 
Recall that
\begin{equation}
\begin{aligned}
m_{ab}(x,y)
&=
\sum_{k=0}^{K-1}\gamma_{abk}(y)\phi_k(x)
+
o(1)K^{-\frac{l_2-1}{2}-\frac{1}{12}} \\
&=
\sum_{k=0}^{K-1}
\left(
\sum_{d=0}^{K-1}\gamma_{abkd}\phi_d(y)
+
o(1)K^{-\frac{l_2-1}{2}-\frac{1}{12}}
\right)\phi_k(x)
+
o(1)K^{-\frac{l_2-1}{2}-\frac{1}{12}} \\
&=
\Phi(x)\otimes\Phi(y)\bm\gamma_{ab}
+
o(1)K^{-\frac{l_2-1}{2}+\frac{11}{12}}.
\end{aligned}
\end{equation}
Therefore, \(\xi_{iab}^2=o(1)K^{-l_2+\frac{17}{6}}\).

Since \(Q<\infty\) is a constant, the summation in \autoref{eqn:Xileq} is finite. 
Consequently,
\begin{equation}
\Xi_i^2
\le
o(1)\left(
K^{-l_1+\frac{5}{6}}
+
K^{-l_2+\frac{17}{6}}
\right)
=
o(1),
\quad i=1,2,\cdots,n,
\end{equation}
as \(l_1\ge 1\), \(l_2\ge 3\), and \(K_n\to\infty\). 
For simplicity, write \(\Xi_i^2\le a_K^2\) uniformly in \(i\), where
\begin{equation}
a_K^2
=
o(1)\left(
K^{-l_1+\frac{5}{6}}
+
K^{-l_2+\frac{17}{6}}
\right).
\end{equation}

We now show that the contribution of \(\bm\Xi\) to the score is negligible under finite-variance projections. 
For any \(\bm\vartheta\in\mathbb R^{P_Z}\) with finite \(\bm\vartheta^{\intercal}\bm V\bm\vartheta\), we have
\begin{equation}
\bm\vartheta^{\intercal}
\left(
\bm W^{(n)}-\widetilde{\bm W}^{(n)}
\right)
=
\frac{1}{\sqrt n}
\sum_{i=1}^n
\bm\vartheta^{\intercal}\bm z_i\Xi_i .
\end{equation}
Under the i.i.d. setting and the column-norm bound in Lemma \autoref{lemZ},
\begin{equation}
\begin{aligned}
\mathrm{Var}
\left\{
\bm\vartheta^{\intercal}
\left(
\bm W^{(n)}-\widetilde{\bm W}^{(n)}
\right)
\right\}
&=
\mathrm{Var}
\left(
\frac{1}{\sqrt n}
\sum_{i=1}^n
\bm\vartheta^{\intercal}\bm z_i\Xi_i
\right) \\
&=
\mathrm{Var}
\left(
\bm\vartheta^{\intercal}\bm z_1\Xi_1
\right) \\
&\le
E\left[
\left(
\bm\vartheta^{\intercal}\bm z_1
\right)^2
\Xi_1^2
\right] \\
&\le
a_K^2
E\left[
\left(
\bm\vartheta^{\intercal}\bm z_1
\right)^2
\right] \\
&=
a_K^2
\bm\vartheta^{\intercal}\bm V\bm\vartheta
=
o(1).
\end{aligned}
\end{equation}
Therefore,
\begin{equation}
\bm\vartheta^{\intercal}
\left(
\bm W^{(n)}-\widetilde{\bm W}^{(n)}
\right)
=
o_P(1).
\end{equation}
Combining this projected negligibility with the projected central limit theorem for \(\widetilde{\bm W}^{(n)}\), we conclude that
\begin{equation}
\bm W^{(n)}
=
\bm W\cdot\{1+o_P(1)\},
\end{equation}
in the sense defined in the lemma.

(b) The vector \(\bm W_1^{(n)}\) is obtained by taking the first \(P_1\) entries of \(\bm W^{(n)}\). 
Therefore, applying part (a) to these entries and splitting the covariance matrix \(\bm V\) into its active and inactive blocks immediately yields
\begin{equation}
\bm W_1^{(n)}
=
\bm W_1\cdot\{1+o_P(1)\},
\end{equation}
with the same projected interpretation stated in part (b) of the lemma.
\qed

\begin{lem}\label{lemPgrow}
Under Assumptions 1-2, and provided $P \le O(n^{c_p})$ for some $c_p > 0$, there exists $k > 0$ satisfying $E(\epsilon_i^{2k}) < \infty$ such that $P_1 = O(n^{c_1})$ and $P_2 = o\left(n^{\left(c_2-c_1\right)k}\right)$.
\end{lem}

\prf
The purpose of this lemma is to verify that the dimension of the expanded linear representation satisfies the growth requirements used in the high-dimensional LASSO argument. 
Condition (a) in Assumption 2 immediately gives $P_1 \le O(n^{c_1})$, since the number of active coefficients is of order $O(K_{n}^2)$ when pairwise interaction components are included. 

For the inactive part, the largest candidate dictionary is generated by the full set of pairwise interactions. 
Thus, the number of inactive coefficients is bounded by
\begin{equation}
P_2 = O(P^2 \cdot K_n^2).
\end{equation}
Since $P \le O(n^{c_p})$ and $K_n^2 \le O(n^{c_1})$, it follows that
\begin{equation}
P_2 = O(P^2 \cdot K_n^2) \le O(n^{2c_p+c_1}).
\end{equation}
Therefore, it is enough to choose $k>0$ such that
\begin{equation}
2c_p+c_1 < (c_2-c_1)k.
\end{equation}
Such a choice is possible because $c_2>c_1$. 
Furthermore, the exponential tail condition in Assumption 1 implies the existence of moments of all finite orders. 
Hence, after choosing $k$ large enough to satisfy the preceding inequality, we still have $E(\epsilon_i^{2k}) < \infty$. 
Consequently, $P_2 = o\left(n^{(c_2-c_1)k}\right)$, and the result follows.
\qed

\begin{lem}\label{lem:event}
For step 1 estimator $\widetilde{\boldsymbol{\beta}}$, we have $\left\{ \widetilde{\boldsymbol{\beta }}=_s\boldsymbol{\beta } \right\} \supset A_n\cap B_n$, where
\begin{equation}
\begin{aligned}
& A_n =
\left\{
\left|  \left(\bm V_{11}^{(n)}\right)^{-1} \bm W^{(n)}_1  \right|
<
\sqrt{n}
\left(
\left|\bm \beta_{\bm P_1}\right|
-
\frac{\lambda_n}{2n}
\left|
\left(\bm V_{11}^{(n)}\right)^{-1}
\mathrm{sign}\left( \bm \beta_{\bm P_1} \right)
\right|
\right)
\right\},
\\
& B_n = \left\{\left|  \left(\bm V_{21}^{(n)}\right) \left(\bm V_{11}^{(n)}\right)^{-1} \bm W^{(n)}_1 -\bm W^{(n)}_2  \right|< \frac{\lambda_n}{2\sqrt{n}}\left(\bm 1 - \left|\left(\bm V_{21}^{(n)}\right)\left(\bm V_{11}^{(n)}\right)^{-1}\mathrm{sign}\left( \bm \beta_{\bm P_1} \right)  \right|  \right)\right\}.
\end{aligned}
\end{equation}
All inequalities are understood componentwise.
\end{lem}

\prf 
Let $\widetilde {\bm u} = \widetilde{ \bm \beta} - \bm \beta$. 
The step 1 LASSO criterion can be written as a function of $\widetilde{\bm u}$:
\begin{equation}
\begin{aligned}
V\left(\widetilde{\bm u}\right) 
& = L\left(\widetilde{\bm \beta} \right)  
= \left\| \bm Y - \bm Z \widetilde{\bm \beta} \right\|_2^2 + \lambda_n  \left\| \widetilde{\bm \beta} \right\|_1 
\\
& = \left\| \bm Z \bm \beta + \bm \Xi + \bm U - \bm Z \widetilde{\bm \beta} \right\|_2^2 + \lambda_n  \left\| \widetilde{\bm \beta} \right\|_1
\\
& = \left\| \bm Z \widetilde{\bm u} - \bm \Xi - \bm U  \right\|_2^2 + \lambda_n  \left\| \widetilde{\bm u} + \bm \beta \right\|_1.
\end{aligned}
\end{equation}
Thus, $\widetilde {\bm u}$ is a minimizer of $V\left({\bm u}\right)$. 
Let $\boldsymbol{\mathcal{V}} = \bm \Xi + \bm U$. 
Then the quadratic part of $V\left( \bm u\right)$ can be rewritten as
\begin{equation}
\left\|\bm Z \bm u - \bm \Xi - \bm U \right\|_2^2 = \left\|\bm Z \bm u - \boldsymbol{\mathcal{V}} \right\|_2^2 
= \bm u^{\intercal} \bm Z^{\intercal} \bm Z \bm u - 2 \bm u^{\intercal} \bm Z^{\intercal} \bm{\mathcal V} + \boldsymbol{\mathcal V}^{\intercal} \boldsymbol{\mathcal V}.
\end{equation}
Since $\bm u^{\intercal} \bm Z^{\intercal} \bm Z \bm u = \left(\sqrt{n}\bm u\right)^{\intercal} \bm V^{(n)} \left(\sqrt{n}\bm u\right)$, we have
\begin{equation}
\begin{aligned}
V\left(\bm u\right) 
= & \left[\left(\sqrt{n}\bm u\right)^{\intercal} \bm V^{(n)} \left(\sqrt{n}\bm u\right) - 2 \bm u^{\intercal} \bm Z^{\intercal} \bm{\mathcal V}  \right]+ \lambda_n  \left\| {\bm u} + \bm \beta \right\|_1+ \boldsymbol{\mathcal V}^{\intercal} \boldsymbol{\mathcal V}
\\
= & V_1\left(\bm u\right) + V_2\left(\bm u\right) +\boldsymbol{\mathcal V}^{\intercal} \boldsymbol{\mathcal V}.
\end{aligned}
\end{equation}
Here, $\boldsymbol{\mathcal V}^{\intercal} \boldsymbol{\mathcal V}$ is independent of $\bm u$. 
For the differentiable part $V_1$, we have
\begin{equation}
\frac{\mathrm{d}V_1}{\mathrm{d}\bm u} = 2\sqrt{n} \left( \bm V^{(n)}\sqrt{n}\bm u - \bm W^{(n)} \right),
\end{equation}
where $\bm W^{(n)}$ is defined in Lemma \autoref{lemA2}. 

We now construct a candidate solution satisfying the KKT conditions with the correct sign pattern. 
Partition $\bm u$ and $\bm W^{(n)}$ conformably with the active and inactive components. 
If the active part satisfies
\begin{equation}
\left\{ \mathrm{sign}\left(\bm \beta_{\bm P_1}\right) \widetilde{\bm u}_1 > - \left|\bm \beta_{\bm P_1}\right|  \right\},
\end{equation}
then $\widetilde{\bm \beta}_{\bm P_1}=\bm \beta_{\bm P_1}+\widetilde{\bm u}_1$ has the same componentwise signs as $\bm \beta_{\bm P_1}$. 
Therefore,
\begin{equation}
\left\{ \mathrm{sign}\left(\bm \beta_{\bm P_1}\right) \widetilde{\bm u}_1 > - \left|\bm \beta_{\bm P_1}\right|  \right\} \subset \left\{ \widetilde{\bm \beta}_{\bm P_1} =_s \bm \beta_{\bm P_1} \right\}.
\end{equation}

Under this sign pattern and with $\widetilde{\bm \beta}_{\bm P_2}=0$, the KKT conditions reduce to
\begin{equation}
\begin{aligned}
& \bm V_{11}^{(n)}\sqrt{n} \widetilde{\bm u}_1  - \bm W^{(n)}_1 = -\frac{\lambda_n}{2\sqrt{n}} \mathrm{sign}\left( \bm \beta_{\bm P_1} \right),
\\
& \left|\widetilde{\bm u}_1\right| < \left|\bm \beta_{\bm P_1}\right|,
\\
& -\frac{\lambda_n}{2\sqrt{n}} \bm 1 \le  \bm V_{21}^{(n)} \sqrt{n} \widetilde{\bm u}_1 -\bm W^{(n)}_2 \le  \frac{\lambda_n}{2\sqrt{n}} \bm 1.
\end{aligned}
\end{equation}
The first equation determines the active perturbation:
\begin{equation}
\sqrt{n}\widetilde{\bm u}_1
=
\left(\bm V_{11}^{(n)}\right)^{-1}\bm W_1^{(n)}
-
\frac{\lambda_n}{2\sqrt n}
\left(\bm V_{11}^{(n)}\right)^{-1}
\mathrm{sign}\left( \bm \beta_{\bm P_1} \right).
\end{equation}
The event $A_n$ guarantees that the resulting active perturbation is small enough, component by component, so that $\left|\widetilde{\bm u}_1\right| < \left|\bm \beta_{\bm P_1}\right|$ and hence the active signs are preserved.

It remains to verify the inactive KKT inequalities. 
Substituting the expression for $\sqrt{n}\widetilde{\bm u}_1$ into the inactive score gives
\begin{equation}
\begin{aligned}
\bm V_{21}^{(n)} \sqrt{n} \widetilde{\bm u}_1 -\bm W^{(n)}_2
&=
\bm V_{21}^{(n)}
\left(\bm V_{11}^{(n)}\right)^{-1}
\bm W_1^{(n)}
-\bm W_2^{(n)}
\\
&\quad
-
\frac{\lambda_n}{2\sqrt n}
\bm V_{21}^{(n)}
\left(\bm V_{11}^{(n)}\right)^{-1}
\mathrm{sign}\left( \bm \beta_{\bm P_1} \right).
\end{aligned}
\end{equation}
By the triangle inequality, the inactive KKT inequalities hold if
\begin{equation}
\left|  \left(\bm V_{21}^{(n)}\right) \left(\bm V_{11}^{(n)}\right)^{-1} \bm W^{(n)}_1 -\bm W^{(n)}_2  \right|
<
\frac{\lambda_n}{2\sqrt{n}}
\left(\bm 1 - \left|\left(\bm V_{21}^{(n)}\right)\left(\bm V_{11}^{(n)}\right)^{-1}\mathrm{sign}\left( \bm \beta_{\bm P_1} \right)  \right|  \right),
\end{equation}
which is precisely the event $B_n$.

Therefore, on $A_n\cap B_n$, there exists a vector satisfying the KKT conditions with $\widetilde{\bm \beta}_{\bm P_1} =_s \bm \beta_{\bm P_1}$ and $\widetilde{\bm \beta}_{\bm P_2}=0$. 
Since the LASSO objective is convex, these KKT conditions are sufficient for optimality. 
Hence,
\begin{equation}
A_n\cap B_n \subset \left\{ \widetilde{\boldsymbol{\beta }}=_s\boldsymbol{\beta } \right\},
\end{equation}
which proves the result.
\qed

\section{Proofs}

\subsection{Proof of {Theorem 4.1}}

\prf 
The proof proceeds by reducing the selection consistency of the first-step LASSO estimator to the sign-consistency result of \cite{zhao2006model}. 
By Lemma \autoref{lem:event}, it is enough to show that the events $A_n$ and $B_n$ occur with probability tending to one. 
We claim that
\begin{equation}
P\left( A_n^c\right) + P\left( B_n^c \right) = O\left(\frac{P_Z n^k}{\lambda_n^{2k}}\right) \to 0.
\end{equation}
Once this claim is established, Lemma \autoref{lem:event} immediately implies
$\widetilde{\bm \beta}_{\bm P_Z} =_s \bm \beta_{\bm P_Z}$ with probability tending to one.

We now verify the conditions required by Theorem 3 in \cite{zhao2006model}. 
In our notation, these conditions are:
\begin{enumerate}
    \item[(i)] $\frac{1}{n}Z_j^{\intercal}Z_j \le M_1$ for all $j$;
    \item[(ii)] $\rho_{\min}(\boldsymbol{V}_{11})>c_3>0$;
    \item[(iii)] $P_1=O(n^{c_1})$;
    \item[(iv)] $n^{\frac{1-c_2}{2}} d(\boldsymbol{\beta}_{P_1})\ge c_4>0$;
    \item[(v)] the strong irrepresentable condition holds;
    \item[(vi)] $P_2=o(n^{(c_2-c_1)k})$ for some $k$ satisfying $E(\epsilon_i^{2k})<\infty$;
    \item[(vii)] $\frac{1}{P_Z}\left(  \frac{\lambda_n}{\sqrt{n}}\right)^{2k} \to \infty$.
\end{enumerate}

We check these conditions one by one. 
Condition (i) follows from Lemma \autoref{lemZ}, which controls the column norms of the expanded design matrix. 
Conditions (ii) and (iv) are exactly imposed in Assumption 2. 
Conditions (iii) and (vi) follow from Lemma \autoref{lemPgrow}, which translates the growth of the original covariate dimension and the truncation level into the growth rates of the expanded active and inactive dimensions. 
Condition (v) is Assumption 3. 
Finally, Assumption 1 gives finite moments of all orders because the error distribution has an exponential-type tail. 
Together with the imposed tuning-rate condition on $\lambda_n$, this verifies condition (vii). 
Equivalently, condition (vii) is the same as $P_Z n^k/\lambda_n^{2k}\to 0$.

Therefore, by Theorem 3 of \cite{zhao2006model},
\begin{equation}
P\left( A_n^c\right) + P\left( B_n^c \right)
=
O\left( \frac{P_Z n^k}{\lambda_n^{2k}} \right)
\to 0.
\end{equation}
Combining this bound with Lemma \autoref{lem:event}, we obtain
\begin{equation}
P\left( \widetilde{\bm \beta}_{\bm P_Z} \neq_s \bm \beta_{\bm P_Z} \right)
\le
P\left(\left(A_n \cap B_n \right)^C\right)
\le
P\left(A_n^C\right)+P\left(B_n^C\right)
=
O\left( \frac{P_Z n^k}{\lambda_n^{2k}} \right)
\to 0.
\end{equation}
This establishes the sign consistency of the first-step estimator. 
Since the second-step adaptive group LASSO assigns infinite weights to components estimated as zero in the first step, the same zero pattern is preserved in Step 2. 
Hence the selected model is sign consistent as stated in {Theorem 4.1}.
\qed

\subsection{Proof of {Theorem 4.2}}

\prf 
(a) We first study the active-set estimator obtained in Step 2. 
On the sign-consistency event from Step 1, the inactive components are removed by the infinite adaptive weights. 
Therefore, Step 2 reduces to an adaptive group LASSO problem on the active set.

We define
\begin{equation}
\begin{cases}
\tau_{\alpha_l}
=
\widehat{w}_{\alpha_l}\mathrm{sign}\left(\widehat{\alpha}_l\right),
& l\in \boldsymbol{L},\\
\tau_{\boldsymbol{b}_r}
=
\widehat{w}_{\boldsymbol{b}_r}
\frac{\widehat{\boldsymbol{b}}_r
\mathbf{1}_{\left\{\left\|\widehat{\boldsymbol{b}}_r\right\|_2>0\right\}}}
{\left\|\widehat{\boldsymbol{b}}_r\right\|_2},
& r\in \boldsymbol{R},\\
\tau_{\boldsymbol{\gamma}_{ab}}
=
\widehat{w}_{\boldsymbol{\gamma}_{ab}}
\frac{\widehat{\boldsymbol{\gamma}}_{ab}
\mathbf{1}_{\left\{\left\|\widehat{\boldsymbol{\gamma}}_{ab}\right\|_2>0\right\}}}
{\left\|\widehat{\boldsymbol{\gamma}}_{ab}\right\|_2},
& (a,b)\in \boldsymbol{S},
\end{cases}
\end{equation}
which are the sub-gradients of the penalty terms for the active set in Step 2.
Under the event $C_n$, we have $\widehat{\boldsymbol{\beta}}_{P_2}=\mathbf{0}$. 
Hence $\widehat{\boldsymbol{\beta}}_{P_1}$ solves the active-set normal equation
\begin{equation}
-2\boldsymbol{Z}_{P_1}^{\intercal}
\left(
\boldsymbol{Y}-\boldsymbol{Z}_{P_1}\widehat{\boldsymbol{\beta}}_{P_1}
\right)
+
\widetilde{\lambda}_n\boldsymbol{\tau}_n
=
0,
\end{equation}
where $\boldsymbol{\tau}_n$ collects all active-set sub-gradients. 
Substituting 
$
\boldsymbol{Y}
=
\boldsymbol{Z}_{P_1}\boldsymbol{\beta}_{P_1}
+
\boldsymbol{U}
+
\boldsymbol{\Xi}
$
into the preceding display gives
\begin{equation}
2\left(\boldsymbol{Z}_{P_1}^{\intercal}\boldsymbol{Z}_{P_1}\right)
\left(\widehat{\boldsymbol{\beta}}_{P_1}-\boldsymbol{\beta}_{P_1}\right)
-
2\boldsymbol{Z}_{P_1}^{\intercal}
\left(\boldsymbol{U}+\boldsymbol{\Xi}\right)
+
\widetilde{\lambda}_n\boldsymbol{\tau}_n
=
0,
\end{equation}
with probability tending to one. 
Rearranging yields
\begin{equation}\label{eqn:prfbetan}
\begin{aligned}
\sqrt{n}
\left(
\widehat{\boldsymbol{\beta}}_{P_1}
-
\boldsymbol{\beta}_{P_1}
\right)
&=
n
\left(
\boldsymbol{Z}_{P_1}^{\intercal}\boldsymbol{Z}_{P_1}
\right)^{-1}
\boldsymbol{W}_{1}^{(n)}
-
\frac{1}{2}
\sqrt{n}
\left(
\boldsymbol{Z}_{P_1}^{\intercal}\boldsymbol{Z}_{P_1}
\right)^{-1}
\widetilde{\lambda}_n
\boldsymbol{\tau}_n
\\
&=
\left(\boldsymbol{V}_{11}^{(n)}\right)^{-1}
\boldsymbol{W}_{1}^{(n)}
-
\frac{1}{2}
\frac{\widetilde{\lambda}_n}{\sqrt{n}}
\left(\boldsymbol{V}_{11}^{(n)}\right)^{-1}
\boldsymbol{\tau}_n,
\end{aligned}
\end{equation}
where
$
\boldsymbol{W}_{1}^{(n)}
=
n^{-1/2}
\boldsymbol{Z}_{P_1}^{\intercal}
\left(
\boldsymbol{U}+\boldsymbol{\Xi}
\right)
$
is defined in Lemma \autoref{lemA2}. 
The first term in \autoref{eqn:prfbetan} is the usual least-squares score on the active set, while the second term is the bias induced by the adaptive group penalty. 
The rest of the proof shows that this bias term is asymptotically negligible.

We first verify that the active Gram matrix is stable. 
By independence and the boundedness of the basis functions,
\begin{equation}
\begin{aligned}
E\left\|
\boldsymbol{V}_{11}^{(n)}
-
\boldsymbol{V}_{11}
\right\|_2^2
&=
\frac{1}{n^2}
E\left\|
\boldsymbol{Z}_{P_1}^{\intercal}\boldsymbol{Z}_{P_1}
-
E\left(
\boldsymbol{Z}_{P_1}^{\intercal}\boldsymbol{Z}_{P_1}
\right)
\right\|_2^2
\\
&=
\frac{1}{n^2}
\sum_{i=1}^n
E\left\|
\boldsymbol{z}_{iP_1}^{\intercal}\boldsymbol{z}_{iP_1}
-
E\left[
\boldsymbol{z}_{iP_1}^{\intercal}\boldsymbol{z}_{iP_1}
\right]
\right\|_2^2
\\
&\le
\frac{CQK^4 n}{n^2}
=
O\left(\frac{K^4}{n}\right).
\end{aligned}
\end{equation}
This term converges to zero when $K_n=o(n^{1/4})$. 
Together with Assumption 2(c), this implies that 
$
\left(\boldsymbol{V}_{11}^{(n)}\right)^{-1}
$
is well behaved with probability tending to one.

It remains to control $\boldsymbol{\tau}_n$. 
The adaptive weights in Step 2 are constructed from the Step 1 estimator. 
Thus, we inspect the Step 1 active-set normal equation. 
Under the event $A_n\cap B_n$, the first-step estimator has the correct signs on the active set, and the KKT condition gives
\begin{equation}
-2\boldsymbol{Z}_{P_1}^{\intercal}
\left(
\boldsymbol{Y}
-
\boldsymbol{Z}_{P_1}
\widetilde{\boldsymbol{\beta}}_{P_1}
\right)
+
\lambda_n
\mathrm{sign}
\left(
\boldsymbol{\beta}_{P_1}
\right)
=
0,
\end{equation}
where the sub-gradient equals 
$
\mathrm{sign}\left(\boldsymbol{\beta}_{P_1}\right)
$
because the coefficients are not grouped in Step 1.
Solving this equation gives
\begin{equation}
\sqrt{n}
\left(
\widetilde{\boldsymbol{\beta}}_{P_1}
-
\boldsymbol{\beta}_{P_1}
\right)
=
\left(
\boldsymbol{V}_{11}^{(n)}
\right)^{-1}
\boldsymbol{W}_{1}^{(n)}
-
\frac{\lambda_n}{2\sqrt{n}}
\left(
\boldsymbol{V}_{11}^{(n)}
\right)^{-1}
\mathrm{sign}
\left(
\boldsymbol{\beta}_{P_1}
\right).
\end{equation}

We next show that the first-step active coefficients are uniformly separated from zero. 
For a vector $\boldsymbol{g}$, by the union bound,
\begin{equation}
\begin{aligned}
P\left(
\left|
\widetilde{\boldsymbol{\beta}}_{P_1}
-
\boldsymbol{\beta}_{P_1}
\right|
>
\boldsymbol{g}
\right)
&\le
\sum_{j=1}^{P_1}
P\left(
\left|
\widetilde{\beta}_j-\beta_j
\right|
>
g_j
\right)
\\
&=
\sum_{j=1}^{P_1}
P\left(
\sqrt{n}
\left|
\widetilde{\beta}_j-\beta_j
\right|
>
\sqrt{n}g_j
\right)
\\
&=
\sum_{j=1}^{P_1}
P\left(
\left|
\zeta_j^{(n)}
-
\frac{\lambda_n}{2\sqrt{n}}\varpi_j
\right|
>
\sqrt{n}g_j
\right)
\\
&\le
\sum_{j=1}^{P_1}
P\left(
\left|
\zeta_j^{(n)}
\right|
>
\sqrt{n}
\left(
g_j+\frac{\lambda_n}{2n}\varpi_j
\right)
\right),
\end{aligned}
\end{equation}
where
$
\boldsymbol{\beta}_{P_1}
=
(\beta_1,\cdots,\beta_{P_1})^{\intercal}
$,
$
\left(\boldsymbol{V}_{11}^{(n)}\right)^{-1}
\boldsymbol{W}_{1}^{(n)}
=
\boldsymbol{\zeta}^{(n)}
$,
and
$
\left(\boldsymbol{V}_{11}^{(n)}\right)^{-1}
\mathrm{sign}\left(\boldsymbol{\beta}_{P_1}\right)
=
(\varpi_1,\cdots,\varpi_{P_1})^{\intercal}
$.
Since $\zeta_j^{(n)}$ has finite $2k$-th order moment, Markov's inequality gives
\begin{equation}
\begin{aligned}
P\left(
\left|
\zeta_j^{(n)}
\right|
>
\sqrt{n}
\left(
g_j+\frac{\lambda_n}{2n}\varpi_j
\right)
\right)
&\le
O\left(
\left[
\sqrt{n}
\left(
g_j+\frac{\lambda_n}{2n}\varpi_j
\right)
\right]^{-2k}
\right)
\\
&=
O\left(
n^{-k}g_j^{-2k}
+
\frac{n^k}{\lambda_n^{2k}}
\right).
\end{aligned}
\end{equation}
Therefore,
\begin{equation}
P\left(
\left|
\widetilde{\boldsymbol{\beta}}_{P_1}
-
\boldsymbol{\beta}_{P_1}
\right|
>
\boldsymbol{g}
\right)
\le
O\left(
K^2
\left[
\frac{\max\{g_j^{-2k}\}}{n^k}
+
\frac{n^k}{\lambda_n^{2k}}
\right]
\right).
\end{equation}
Choose
$
g_1=g_2=\cdots=g_{P_1}
=
\frac{c_4}{2}
n^{-(1-c_2)/2}
\le
\frac{1}{2}
d\left(\boldsymbol{\beta}_{P_1}\right)
$.
Then
\begin{equation}
P\left(
\left|
\widetilde{\boldsymbol{\beta}}_{P_1}
-
\boldsymbol{\beta}_{P_1}
\right|
>
\boldsymbol{g}
\right)
\le
O\left(
K^2 n^{-kc_2}
\right)
+
O\left(
K^2\frac{n^k}{\lambda_n^{2k}}
\right)
=
o(1),
\end{equation}
where the last equality follows by choosing $k$ sufficiently large and by the tuning-rate condition on $\lambda_n$. 
Thus,
$
d\left(\widetilde{\boldsymbol{\beta}}_{P_1}\right)
\ge
(c_4/2)n^{-(1-c_2)/2}
$
with probability tending to one.

We now return to the adaptive penalty in Step 2. 
Since all active first-step estimates are bounded away from zero at the preceding rate, the adaptive weights on the active set are at most of order $O_P(n^{(1-c_2)/2})$. 
Therefore, componentwise on the active set,
\begin{equation}
\frac{\widetilde{\lambda}_n}{\sqrt{n}}
\left(
\boldsymbol{V}_{11}^{(n)}
\right)^{-1}
\boldsymbol{\tau}_n
=
\mathbf{1}_{P_1}
\cdot
O_P\left(
\frac{\widetilde{\lambda}_n}{\sqrt{n}}
n^{(1-c_2)/2}
\right)
=
\mathbf{1}_{P_1}
\cdot
O_P\left(
n^{-c_2/2}
\widetilde{\lambda}_n
\right)
=
\mathbf{1}_{P_1}
\cdot
o_P(1),
\end{equation}
where the last equality follows from
$
\widetilde{\lambda}_n=o(n^{c_2/2})
$.
Substituting this into \autoref{eqn:prfbetan}, we see that the second-step estimator has the same first-order stochastic order as the active-set least-squares estimator. 
In particular,
\begin{equation}\label{eqn:betadis}
\widehat{\boldsymbol{\alpha}}_{\boldsymbol{L}}
-
\boldsymbol{\alpha}_{\boldsymbol{L}}
=
O_P(n^{-1/2}),
\quad
\left\|
\widehat{\boldsymbol{b}}_r
-
\boldsymbol{b}_r
\right\|_2^2
=
O_P\left(\frac{K}{n}\right),
\quad
\left\|
\widehat{\boldsymbol{\gamma}}_{ab}
-
\boldsymbol{\gamma}_{ab}
\right\|_2^2
=
O_P\left(\frac{K^2}{n}\right).
\end{equation}
Taking the first $L$ entries of $\widehat{\boldsymbol{\beta}}_{P_1}$, we obtain
$
\widehat{\boldsymbol{\alpha}}_{\boldsymbol{L}}
-
\boldsymbol{\alpha}_{\boldsymbol{L}}
=
O_P(n^{-1/2})
=
o_P(1)
$.

(b) We next prove the $L^2$-consistency of the univariate nonparametric components. 
The estimator is
\begin{equation}
\widehat{m}_{r}(x)
=
\Phi(x)\widehat{\boldsymbol{b}}_r
=
\sum_{k=0}^{K-1}
\widehat{b}_{rk}\phi_k(x).
\end{equation}
By orthonormality of the Hermite basis, the squared $L^2$-error decomposes into the estimation error of the first $K$ coefficients and the truncation error:
\begin{equation}
\begin{aligned}
\left\|
m_r-\widehat{m}_r
\right\|_{L^2}^2
&=
\sum_{k=0}^{K-1}
\left(
b_{rk}-\widehat{b}_{rk}
\right)^2
+
\sum_{k\ge K}
b_{rk}^2
\\
&\le
\sum_{k=0}^{K-1}
\left(
b_{rk}-\widehat{b}_{rk}
\right)^2
+
o(1)K^{-l_1},
\end{aligned}
\end{equation}
where the last inequality follows from Lemma \autoref{lemA1}. 
By \autoref{eqn:betadis},
\begin{equation}
\sum_{k=0}^{K-1}
\left(
b_{rk}-\widehat{b}_{rk}
\right)^2
=
O_P\left(\frac{K}{n}\right).
\end{equation}
Hence,
\begin{equation}
\left\|
m_r-\widehat{m}_r
\right\|_{L^2}^2
\le
O_P\left(
\frac{K}{n}
+
K^{-l_1}
\right).
\end{equation}
Since $l_1\ge 1$ and $K=o(n^{1/4})$, both terms vanish, and thus
\begin{equation}
\left\|
m_r-\widehat{m}_r
\right\|_{L^2}^2
=
o_P(1).
\end{equation}
This proves the $L^2$-consistency of $\widehat{m}_r$.

(c) The proof for the interactive components is analogous. 
The tensor-product Hermite basis gives an orthogonal decomposition of the bivariate $L^2$-error into a coefficient-estimation part and a truncation part. 
Using \autoref{eqn:betadis} and Lemma \autoref{lemA1}, we obtain
\begin{equation}
\left\|
m_{ab}
-
\widehat{m}_{ab}
\right\|_{L^2}^2
\le
O_P\left(\frac{K^2}{n}\right)
+
o(1)K^{-l_2+1}.
\end{equation}
The first term vanishes because $K=o(n^{1/4})$ implies $K^2/n=o(1)$. 
The second term vanishes under the smoothness condition $l_2\ge 3$. 
Therefore,
$
\left\|
m_{ab}
-
\widehat{m}_{ab}
\right\|_{L^2}^2=o_P(1)
$.

(d) Finally, setting $K=O(n^\nu)$ in the preceding bounds gives the explicit rates. 
For the univariate components,
\begin{equation}
\left\|
m_r-\widehat{m}_r
\right\|_{L^2}^2
\le
O_P\left(
n^{-(1-\nu)}
+
n^{-l_1\nu}
\right).
\end{equation}
For the interactive components,
\begin{equation}
\left\|
m_{ab}-\widehat{m}_{ab}
\right\|_{L^2}^2
\le
O_P\left(
n^{-(1-2\nu)}
+
n^{-l_2\nu+\nu}
\right).
\end{equation}
Taking square roots of the two displays yields the rates stated in the theorem.
\qed

\subsection{Proof of Theorem 4.3}

\prf 
We start from the expansion in \autoref{eqn:prfbetan}:
\begin{equation}
\sqrt{n}
\left(
\widehat{\boldsymbol{\beta}}_{P_1}
-
\boldsymbol{\beta}_{P_1}
\right)
=
\left(
\boldsymbol{V}_{11}^{(n)}
\right)^{-1}
\boldsymbol{W}_{1}^{(n)}
-
\frac{1}{2}
\frac{\widetilde{\lambda}_n}{\sqrt{n}}
\left(
\boldsymbol{V}_{11}^{(n)}
\right)^{-1}
\boldsymbol{\tau}_n.
\end{equation}
The proof applies the central limit theorem to the fixed-dimensional projected statistic appearing in the theorem, rather than to the whole growing-dimensional coefficient vector. 
This is important because the number of basis coefficients may increase with $n$, whereas the number of reported quantities, namely the linear coefficients and the function values at a fixed point, remains fixed.

We first relate the function-value errors to the coefficient errors. 
For a specific $r_0$, we have
\begin{equation}
\sqrt{n}
\left(
\widehat{m}_{r_0}\left(x_{r_0}\right)
-
m_{r_0}\left(x_{r_0}\right)
\right)
=
\Phi\left(x_{r_0}\right)
\sqrt{n}
\left(
\widehat{\boldsymbol{b}}_{r_0}
-
\boldsymbol{b}_{r_0}
\right)
+
\sqrt{n}\xi_{r_0},
\end{equation}
where $\xi_{r_0}$ is the sieve approximation remainder at $x_{r_0}$. 
By Lemma \autoref{lemA1}, $\xi_{r_0}=o(1)K^{-\frac{l_1}{2}+\frac{5}{12}}$. 
Since $K=O(n^\nu)$ and $\nu \ge 6/(6l_1-5)$, we obtain
\begin{equation}
\sqrt{n}\xi_{r_0}
=
o(1)n^{\frac{1}{2}}
n^{-\frac{l_1\nu}{2}+\frac{5\nu}{12}}
=
o(1)n^{-\frac{(6l_1-5)\nu-6}{12}}
=
o(1).
\end{equation}
Thus, after multiplication by $\sqrt{n}$, the univariate truncation error is asymptotically negligible.

Similarly, for a specific $(a_0,b_0)$, we have
\begin{equation}
\begin{aligned}
&
\sqrt{n}
\left(
\widehat{m}_{a_0b_0}
\left(
x_{a_0},x_{b_0}
\right)
-
m_{a_0b_0}
\left(
x_{a_0},x_{b_0}
\right)
\right)
\\
&\qquad
=
\left(
\Phi\left(x_{a_0}\right)
\otimes
\Phi\left(x_{b_0}\right)
\right)
\sqrt{n}
\left(
\widehat{\boldsymbol{\gamma}}_{a_0b_0}
-
\boldsymbol{\gamma}_{a_0b_0}
\right)
+
\sqrt{n}\xi_{a_0b_0},
\end{aligned}
\end{equation}
where $\xi_{a_0b_0}$ is the bivariate sieve approximation remainder and
$\xi_{a_0b_0}=o(1)K^{-\frac{l_2}{2}+\frac{17}{12}}$. 
Since $K=O(n^\nu)$ and $\nu \ge 6/(6l_2-17)$,
\begin{equation}
\sqrt{n}\xi_{a_0b_0}
=
o(1)n^{\frac{1}{2}}
n^{-\frac{l_2\nu}{2}+\frac{17\nu}{12}}
=
o(1)n^{-\frac{(6l_2-17)\nu-6}{12}}
=
o(1).
\end{equation}
Hence the bivariate truncation error is also negligible at the $\sqrt{n}$ scale.

Collecting the preceding displays for all selected univariate and interactive components gives
\begin{equation}
\begin{aligned}
&
\sqrt{n}
\left(
\widehat{m}_{r_j}
\left(
x_{r_j}
\right)
-
m_{r_j}
\left(
x_{r_j}
\right)
\right)_{1\le j\le R}
\\
&\qquad
=
\mathrm{diag}
\left\{
\left(
\Phi
\left(
x_{r_j}
\right)
\right)_{1\le j\le R}
\right\}
\sqrt{n}
\left(
\widehat{\boldsymbol{b}}_{\boldsymbol{R}}
-
\boldsymbol{b}_{\boldsymbol{R}}
\right)
+
o_P(1),
\\
&
\sqrt{n}
\left(
\widehat{m}_{a_jb_j}
\left(
x_{a_j},x_{b_j}
\right)
-
m_{a_jb_j}
\left(
x_{a_j},x_{b_j}
\right)
\right)_{1\le j\le S}
\\
&\qquad
=
\mathrm{diag}
\left\{
\left(
\Phi
\left(
x_{a_j}
\right)
\otimes
\Phi
\left(
x_{b_j}
\right)
\right)_{1\le j\le S}
\right\}
\sqrt{n}
\left(
\widehat{\boldsymbol{\gamma}}_{\boldsymbol{S}}
-
\boldsymbol{\gamma}_{\boldsymbol{S}}
\right)
+
o_P(1).
\end{aligned}
\end{equation}

Therefore,
\begin{equation}
\sqrt{n}
\left(
\begin{array}{c}
\widehat{\boldsymbol{\alpha}}_{\boldsymbol{L}}
-
\boldsymbol{\alpha}_{\boldsymbol{L}}
\\
\left(
\widehat{m}_{r_j}
\left(
x_{r_j}
\right)
-
m_{r_j}
\left(
x_{r_j}
\right)
\right)_{1\le j\le R}
\\
\left(
\widehat{m}_{a_jb_j}
\left(
x_{a_j},x_{b_j}
\right)
-
m_{a_jb_j}
\left(
x_{a_j},x_{b_j}
\right)
\right)_{1\le j\le S}
\end{array}
\right)
=
\boldsymbol{G}(\bm{x})
\sqrt{n}
\left(
\widehat{\boldsymbol{\beta}}_{P_1}
-
\boldsymbol{\beta}_{P_1}
\right)
+
o_P(1),
\end{equation}
where
\begin{equation}
\boldsymbol{G}(\bm{x})
=
\left(
\begin{matrix}
I_L & O & O\\
O &
\mathrm{diag}
\left\{
\left(
\Phi
\left(
x_{r_j}
\right)
\right)_{1\le j\le R}
\right\}
& O\\
O & O &
\mathrm{diag}
\left\{
\left(
\Phi
\left(
x_{a_j}
\right)
\otimes
\Phi
\left(
x_{b_j}
\right)
\right)_{1\le j\le S}
\right\}
\end{matrix}
\right).
\end{equation}
The matrix $\boldsymbol{G}(\bm{x})$ only evaluates the selected basis expansions at the fixed point $\bm{x}$, so its number of rows is $L+R+S$, which is fixed.

Substituting \autoref{eqn:prfbetan} into the last display gives
\begin{equation}
\begin{aligned}
&
\sqrt{n}
\left(
\begin{array}{c}
\widehat{\boldsymbol{\alpha}}_{\boldsymbol{L}}
-
\boldsymbol{\alpha}_{\boldsymbol{L}}
\\
\left(
\widehat{m}_{r_j}
\left(
x_{r_j}
\right)
-
m_{r_j}
\left(
x_{r_j}
\right)
\right)_{1\le j\le R}
\\
\left(
\widehat{m}_{a_jb_j}
\left(
x_{a_j},x_{b_j}
\right)
-
m_{a_jb_j}
\left(
x_{a_j},x_{b_j}
\right)
\right)_{1\le j\le S}
\end{array}
\right)
\\
&\qquad
=
\boldsymbol{G}(\bm{x})
\left(
\boldsymbol{V}_{11}^{(n)}
\right)^{-1}
\boldsymbol{W}_{1}^{(n)}
-
\frac{1}{2}
\frac{\widetilde{\lambda}_n}{\sqrt{n}}
\boldsymbol{G}(\bm{x})
\left(
\boldsymbol{V}_{11}^{(n)}
\right)^{-1}
\boldsymbol{\tau}_n
+
o_P(1).
\end{aligned}
\end{equation}
We now study the two terms on the right-hand side separately.

First, since $L+R+S$ is fixed, the central limit theorem applies directly to the fixed-dimensional projected score
$
\boldsymbol{G}(\bm{x})
\left(
\boldsymbol{V}_{11}^{(n)}
\right)^{-1}
\boldsymbol{W}_{1}^{(n)}
$.
Lemma \autoref{lemA2} shows that the approximation error in $\boldsymbol{W}_{1}^{(n)}$ is asymptotically negligible, and the preceding proof of Theorem 4.2 gives
$
\boldsymbol{V}_{11}^{(n)}
\overset{P}{\longrightarrow}
\boldsymbol{V}_{11}
$.
Therefore,
\begin{equation}
\left(
\sigma^2
\boldsymbol{G}(\bm{x})
\boldsymbol{V}_{11}^{-1}
\boldsymbol{G}(\bm{x})^{\intercal}
\right)^{-1/2}
\boldsymbol{G}(\bm{x})
\left(
\boldsymbol{V}_{11}^{(n)}
\right)^{-1}
\boldsymbol{W}_{1}^{(n)}
\overset{d}{\longrightarrow}
N
\left(
\boldsymbol{0},
\boldsymbol{I}_{L+R+S}
\right).
\end{equation}

Second, we show that the adaptive penalty term is negligible under the same normalization. 
From the proof of Theorem 4.2, on the active set the sub-gradient vector satisfies
$
\|\boldsymbol{\tau}_n\|_2
=
O_P\left(n^{(1-c_2)/2}\right)
$.
Since
$
\widetilde{\lambda}_n=o(n^{c_2/2})
$,
we have
$
\widetilde{\lambda}_n n^{-1/2} \|\boldsymbol{\tau}_n\|_2=o_P(1)
$.
The matrices
$
\boldsymbol{G}(\bm{x})
$
and
$
\left(
\boldsymbol{V}_{11}^{(n)}
\right)^{-1}
$
are well behaved with probability tending to one, and the normalization matrix is fixed-dimensional. 
Hence,
\begin{equation}
\left(
\sigma^2
\boldsymbol{G}(\bm{x})
\boldsymbol{V}_{11}^{-1}
\boldsymbol{G}(\bm{x})^{\intercal}
\right)^{-1/2}
\frac{\widetilde{\lambda}_n}{\sqrt{n}}
\boldsymbol{G}(\bm{x})
\left(
\boldsymbol{V}_{11}^{(n)}
\right)^{-1}
\boldsymbol{\tau}_n
=
o_P(1).
\end{equation}

Combining the preceding displays and applying Slutsky's theorem yields
\begin{equation}
\left(
\sigma^2
\boldsymbol{G}(\bm{x})
\boldsymbol{V}_{11}^{-1}
\boldsymbol{G}(\bm{x})^{\intercal}
\right)^{-1/2}
\sqrt{n}
\left(
\begin{array}{c}
\widehat{\boldsymbol{\alpha}}_{\boldsymbol{L}}
-
\boldsymbol{\alpha}_{\boldsymbol{L}}
\\
\left(
\widehat{m}_{r_j}
\left(
x_{r_j}
\right)
-
m_{r_j}
\left(
x_{r_j}
\right)
\right)_{1\le j\le R}
\\
\left(
\widehat{m}_{a_jb_j}
\left(
x_{a_j},x_{b_j}
\right)
-
m_{a_jb_j}
\left(
x_{a_j},x_{b_j}
\right)
\right)_{1\le j\le S}
\end{array}
\right)
\overset{d}{\longrightarrow}
N
\left(
\boldsymbol{0},
\boldsymbol{I}_{L+R+S}
\right),
\end{equation}
which proves the result.
\qed

\section{Subgradient Descent Algorithm for S-LASSO}
\label{subgde}

We first find the classical LASSO solution by minimizing
\begin{equation}
\begin{aligned}
\min_{\bm{\alpha}, \bm{b}, \bm{\gamma}} ||\bm{Y}-\sum_{l=0}^{P}\alpha_{l}\bm{X}_{l}-&\sum_{r=1}^{P}{\Phi} (\bm{X}_{r})\bm{b}_{r}-\sum_{a=1}^{P-1}\sum_{b>a}^{P}{\Phi} (\bm{X}_{a})\otimes {\Phi}(\bm{X}_{b})\bm{\gamma} _{ab}||^{2}_{2}\\&+\lambda _n\left( \sum_{l=0}^{P}|\alpha _l|+\sum_{r=1}^{P}|\bm{b} _{r}| +\sum_{a=1}^{P-1}\sum_{b>a}^{P}|\bm{\gamma}_{ab}|\right).
\end{aligned}
\end{equation}
Then we construct the weights by
\begin{equation}
\begin{aligned}
&\widehat{w}_{\alpha _l}=\begin{cases}
	\sqrt{N_L}\left| \widetilde{\alpha }_l \right|^{-1},&		\left| \widetilde{\alpha }_l \right|>0,\\
	\infty ,&		\left| \widetilde{\alpha }_l \right|=0,\\
\end{cases}
\quad \widehat{w}_{\boldsymbol{\beta }_r}=\begin{cases}
	\sqrt{N_R}|| \widetilde{\boldsymbol{b}}_r || _{2}^{-1},&		|| \widetilde{\boldsymbol{b}}_r || _2>0,\\
	\infty ,&		|| \widetilde{\boldsymbol{b}}_r || _2=0,\\
\end{cases}
\\
&\widehat{w}_{\boldsymbol{\gamma }_{ab}}=\begin{cases}
	\sqrt{N_S}\left\| \widetilde{\boldsymbol{\gamma }}_{ab} \right\| _{2}^{-1},&		\left\| \widetilde{\boldsymbol{\gamma }}_{ab} \right\| _2>0,\\
	\infty ,&		\left\| \widetilde{\boldsymbol{\gamma }}_{ab} \right\| _2=0,\\
\end{cases}
\end{aligned}
\end{equation}
where $\widetilde{\alpha}$, $\widetilde{\bm b}$ and $\widetilde{\bm \gamma}$ are estimators in step 1, $N_L=L$, $N_R=R\times K$ and $N_S=S\times K^2$ are the number of coefficients within each group as our group sizes are significantly different. We use group cardinality to control the strength of the penalty.

Then the adaptive group LASSO objective function is 
\begin{equation}
\begin{aligned}
L\left(\bm{\alpha}, \bm{b}, \bm{\gamma}\right) = ||\bm{Y}-\sum_{l=0}^{P}\alpha_{l}\bm{X}_{l}-&\sum_{r=1}^{P}{\Phi} (\bm{X}_{r})\bm{b}_{r}-\sum_{a=1}^{P-1}\sum_{b>a}^{P}{\Phi} (\bm{X}_{a})\otimes {\Phi}(\bm{X}_{b})\bm{\gamma} _{ab}||^{2}_{2}\\&+\widetilde{\lambda}_n\left( \sum_{l=0}^{P}\widehat{w}_{\alpha_l}|\alpha _l|+\sum_{r=1}^{P}\widehat{w}_{\bm{b} _{r}}\left\|\bm{b} _{r}\right\|_2 +\sum_{a=1}^{P-1}\sum_{b>a}^{P}\widehat{w}_{\bm{\gamma}_{ab}}\left\|\bm{\gamma}_{ab}\right\|_2\right).
\end{aligned}
\end{equation}
One may write the traditional OLS part as $\bm Y- \bm Z \bm \beta$, then
\begin{equation}\label{eqn:lossfun}
\begin{aligned}
L\left(\bm \beta\right) = & ||\bm Y- \bm Z \bm \beta||^{2}_{2}+\widetilde{\lambda}_n\left( \sum_{l=0}^{P}\widehat{w}_{\alpha_l}|\alpha _l|+\sum_{r=1}^{P}\widehat{w}_{\bm{b} _{r}}\left\|\bm{b} _{r}\right\|_2 +\sum_{a=1}^{P-1}\sum_{b>a}^{P}\widehat{w}_{\bm{\gamma}_{ab}}\left\|\bm{\gamma}_{ab}\right\|_2\right)
\\
= & L_\mathrm{OLS}\left(\bm \beta \right) + L_\mathrm{Pen}\left(\bm \beta \right).
\end{aligned}
\end{equation}

The gradient of $L_\mathrm{OLS}$ is $-2\boldsymbol{Z}^{\top}\left( \boldsymbol{Y}-\boldsymbol{Z\beta } \right)$ while the subgradient of $L_\mathrm{Pen}$ is given by
\begin{equation}
\begin{aligned}
&\frac{\partial L_{\mathrm{Pen}}}{\partial \alpha _l}=\widetilde{\lambda }_n\widehat{w}_{\alpha _l}\mathrm{sign}\left( \alpha _l \right) ,\quad l=0,1\cdots ,P,
\\
&\frac{\partial L_{\mathrm{Pen}}}{\partial \boldsymbol{b}_r}=\widetilde{\lambda }_n\widehat{w}_{\boldsymbol{b}_r}\frac{\boldsymbol{b}_r\mathbf{1}_{\left\{ \left\| \boldsymbol{b}_r \right\| _2>0 \right\}}}{\left\| \boldsymbol{b}_r \right\| _2},\quad r=1,\cdots ,P,
\\
&\frac{\partial L_{\mathrm{Pen}}}{\partial \boldsymbol{\gamma }_{ab}}=\widetilde{\lambda }_n\widehat{w}_{\boldsymbol{\gamma }_r}\frac{\boldsymbol{\gamma }_{ab}\mathbf{1}_{\left\{ \left\| \boldsymbol{\gamma }_{ab} \right\| _2>0 \right\}}}{\left\| \boldsymbol{\gamma }_{ab} \right\| _2},\quad 1\le a<b\le P,
\end{aligned}
\end{equation}
where we define $\frac{0}{0}=0$.
Actually, the form $\frac{\partial L_{\mathrm{Pen}}}{\partial \alpha _l}=\widetilde{\lambda }_n\widehat{w}_{\alpha _l}\mathrm{sign}\left( \alpha _l \right) = \widetilde{\lambda }_n\widehat{w}_{\alpha _l}\frac{\alpha_l \mathbf{1}_{\left\{ | \alpha_l|>0\right\}}}{|\alpha_l|}$ is consistent with $\boldsymbol{b}_r$ and $\boldsymbol{\gamma }_{ab}$.
Then we combine them together and get the subgradient of $L$, i.e., 
\begin{equation}
\mathrm{subg}\left(\bm \beta\right)=\frac{\partial L}{\partial \bm \beta} = -2 \bm Z^{\top} \left( \bm Y - \bm Z \bm \beta\right) + \frac{\partial L_{\mathrm{Pen}}}{\partial \bm \beta}.
\end{equation}

Provided the initial value $\bm \beta^{(0)}$ (one may use the step 1 estimates), we iterate the coefficients by
\begin{equation}
\boldsymbol{\beta }^{\left( k \right)}=\boldsymbol{\beta }^{\left( k-1 \right)}-\eta_l \cdot \mathrm{subg}\left( \boldsymbol{\beta }^{\left( k-1 \right)} \right) ,\quad k\ge 1,
\end{equation}
for a given learning rate $\eta_l$.
When $\left\| \boldsymbol{\beta }^{\left( k \right)}-\boldsymbol{\beta }^{\left( k-1 \right)} \right\| _2$ is smaller than a given tolerance level, one can get the estimate $\widehat{\bm \beta} = \boldsymbol{\beta }^{\left( k \right)}$.

\newpage
\section{Additional Tables, Figures and Results}
\label{addtabs}

\subsection{Firm Characteristics}

\begin{scriptsize}
\setstretch{1}
\setlength{\tabcolsep}{4pt}
\begin{longtable}{>{\itshape}p{2.2cm}p{5.5cm}>{\itshape}p{2.2cm}p{5.5cm}}
\caption{\setstretch{1}
Listing of Firm Characteristics} \label{tab:var} \\ 
\toprule
{\rm{\textbf{Acronym}}} & {\textbf{Firm characteristic}} & {\rm{\textbf{Acronym}}} & {\textbf{Firm characteristic}}    \\
\midrule
\endfirsthead

\toprule
{\rm{\textbf{Acronym}}} & {\textbf{Firm characteristic}} & {\rm{\textbf{Acronym}}} & {\textbf{Firm characteristic}}    \\
\midrule
\endhead

\bottomrule
\endfoot

absacc                               & Absolute accruals                                & operprof                             & Operating profitability                             \\
acc                                  & Working capital accruals                         & orgcap                               & Organizational capital                              \\
aeavol                               & Abnormal earnings announcement volume            & pchcurrat                            & \% change in current ratio                          \\
age                                  & \# years since first Compustat coverage          & pchdepr                              & \% change in depreciation                           \\
agr                                  & Asset growth                                     & pchgm\_pchsale                       & \% change in gross margin - \% change in sales      \\
baspread                             & Bid-ask spread                                   & pchquick                             & \% change in quick ratio                            \\
beta                                 & Beta                                             & pchsale\_pchinvt                     & \% change in sales - \% change in inventory         \\
betasq                               & Beta squared                                     & pchsale\_pchrect                     & \% change in sales - \% change in A/R               \\
bm                                   & Book-to-market                                   & pchsale\_pchxsga                     & \% change in sales - \% change in SG\&A             \\
cash                                 & Cash holdings                                    & pchsaleinv                           & \% change sales-to-inventory                        \\
cashdebt                             & Cash flow to debt                                & pctacc                               & Percent accruals                                    \\
cashpr                               & Cash productivity                                & pricedelay                           & Price delay                                         \\
cfp                                  & Cash-flow-to-price ratio                         & ps                                   & Financial statements score                          \\
chcsho                               & Change in shares outstanding                     & quick                                & Quick ratio                                         \\
chfeps                               & Change in forecasted EPS                         & rd                                   & R\&D increase                                       \\
chinv                                & Change in inventory                              & rd\_mve                              & R\&D to market capitalization                       \\
chmom                                & Change in 6-month momentum                       & rd\_sale                             & R\&D to sales                                       \\
chnanalyst                           & Change in number of analysts                     & realestate                           & Real estate holdings                                \\
chtx                                 & Change in tax expense                            & retvol                               & Return volatility                                   \\
cinvest                              & Corporate investment                             & roaq                                 & Return on assets                                    \\
convind                              & Convertible debt indicator                       & roavol                               & Earnings volatility                                 \\
currat                               & Current ratio                                    & roeq                                 & Return on equity                                    \\
depr                                 & Depreciation / PP\&E                             & roic                                 & Return on invested capital                          \\
disp                                 & Dispersion in forecasted EPS                     & rsup                                 & Revenue surprise                                    \\
divi                                 & Dividend initiation                              & salecash                             & Sales to cash                                       \\
divo                                 & Dividend omission                                & saleinv                              & Sales to inventory                                  \\
dolvol                               & Dollar trading volume                            & salerec                              & Sales to receivables                                \\
dy                                   & Dividend to price                                & secured                              & Secured debt                                        \\
ear                                  & Earnings announcement return                     & securedind                           & Secured debt indicator                              \\
egr                                  & Growth in common shareholder equity              & sfe                                  & Scaled earnings forecast                            \\
ep                                   & Earnings to price                                & sgr                                  & Sales growth                                        \\
fgr5yr                               & Forecasted growth in 5-year EPS                  & sin                                  & Sin stocks                                          \\
gma                                  & Gross profitability                              & sp                                   & Sales to price                                      \\
grcapx                               & Growth in capital expenditures                   & std\_dolvol                          & Volatility of liquidity (dollar trading volume)     \\
grltnoa                              & Growth in long-term net operating assets         & std\_turn                            & Volatility of liquidity (share turnover)            \\
herf                                 & Industry sales concentration                     & stdacc                               & Accrual volatility                                  \\
hire                                 & Employee growth rate                             & stdcf                                & Cash flow volatility                                \\
idiovol                              & Idiosyncratic return volatility                  & sue                                  & Unexpected quarterly earnings                       \\
ill                                  & Illiquidity                                      & tang                                 & Debt capacity/firm tangibility                      \\
invest                               & Capital expenditures and inventory               & tb                                   & Tax income to book income                           \\
ipo                                  & New equity issue                                 & turn                                 & Share turnover                                      \\
lev                                  & Leverage                                         & zerotrade                            & Zero trading days                                   \\
lgr                                  & Growth in long-term debt                         & bm\_ia                               & Industry-adjusted book to market                    \\
maxret                               & Maximum daily return                             & cfp\_ia                              & Industry-adjusted cash-flow-to-price ratio          \\
mom12m                               & 12-month momentum                                & chatoia                              & Industry-adjusted change in asset turnover          \\
mom1m                                & 1-month momentum                                 & chempia                              & Industry-adjusted change in employees               \\
mom36m                               & 36-month momentum                                & chpmia                               & Industry-adjusted change in profit margin           \\
mom6m                                & 6-month momentum                                 & pchcapx\_ia                          & Industry adjusted \% change in expenditures \\
ms                                   & Financial statement score                        & indmom                               & Industry momentum                                   \\
mve                                  & Size                                             & mve\_ia                              & Industry-adjusted size                              \\
nanalyst                             & Number of analysts covering stock                & nincr                              &   Number of earnings increases                                                      \\ 

\end{longtable}
\end{scriptsize}
\newpage

\subsection{Empirical Results for Full Specification}

In the main empirical analysis, we restrict the prior interaction set to pairs involving firm size (`mve'). 
This choice is made for two reasons. 
First, it is computationally motivated. With $P=102$ characteristics, the full pairwise interaction set contains $102\times 101/2 = 5{,}151$ candidate interaction terms even before accounting for the linear and nonlinear components. 
In the rolling-window implementation, the resulting design matrix and, in particular, the associated Gram matrix become prohibitively large in practice.
Second, it is economically motivated. In the asset pricing literature, firm size is one of the most extensively studied characteristics, and its interactions with other characteristics are widely documented \citep{HongLimStein2000JF}. 
This makes `mve' a natural and economically meaningful variable to prioritize when introducing interaction terms.

To assess whether this prior-set implementation omits important interaction patterns, we conduct a supplementary full-specification check. 
Specifically, we split the full sample period from 1980 to 2019 into ten consecutive four-year subsamples, and in each subsample we estimate S-LASSO with the complete set of pairwise interactions. 
For each characteristic, we then record how often it is identified as participating in at least one selected interaction effect across the ten subsamples. \autoref{tab:all_cross} summarizes the characteristics most frequently selected in this exercise.

\begin{table}[h!]
\centering
\caption{\setstretch{1.2}
Summary of Characteristics Selected for Interaction Effects.
\footnotesize{This table reports the characteristics most frequently selected as having interaction effects across the ten sub-samples from 1980--1983 to 2016--2019. The reported frequency indicates the number of sub-samples in which the corresponding characteristic is identified as having at least one interaction effect.}}
\label{tab:all_cross}
\footnotesize
\tabcolsep=1.2cm
\setstretch{1.2}
\begin{tabular}{lll}
\toprule
\textbf{\textit{\#Selected}} & \textbf{} & \textbf{Characteristics} \\
\midrule
\textit{\#Selected} = 10 & & mve \\
\textit{\#Selected} = 9  & & indmom \\
\textit{\#Selected} = 6  & & baspread, mom1m, chmom \\
\textit{\#Selected} = 5  & & maxret, turn, std\_turn, dolvol \\
\textit{\#Selected} = 4  & & std\_dolvol, retvol, mom12m, idiovol, zerotrade \\
\bottomrule
\end{tabular}
\end{table}

Several conclusions emerge from the table. 
First, `mve' is selected in all ten subsamples, making it the most frequently connected characteristic in the full specification. 
Second, `indmom' is selected in nine subsamples, and a small group of variables, including `baspread', `mom1m', and `chmom', are also selected repeatedly. 
Third, beyond these leading variables, the frequencies decline quickly, indicating that strong interaction patterns are concentrated in a relatively small subset of characteristics rather than being spread broadly across the full 102-dimensional characteristic space.

Overall, these results provide ex post support for the prior interaction set used in the main text. 
They show that firm size is not only economically well motivated, but also empirically the most central variable in the full-specification interaction network. 
Therefore, focusing on interactions with `mve' in the main analysis captures the most persistent interaction structure while keeping the rolling-window implementation computationally feasible.

\newpage
\subsection{Discussion of the Irrepresentable Condition}
\label{app:irrepresentable}

Assumption~3 is imposed to establish global sign consistency in Theorem~4.1, which may not hold in the case the covariates are highly correlated.
The purpose of this part is therefore not to claim that the condition must literally hold in the empirical application, but to assess how restrictive it appears under dependence structures calibrated to our data.

\subsubsection{Empirical Distribution of Pairwise Correlations}

\begin{figure}[htbp]
    \centering
    \includegraphics[width=0.72\linewidth]{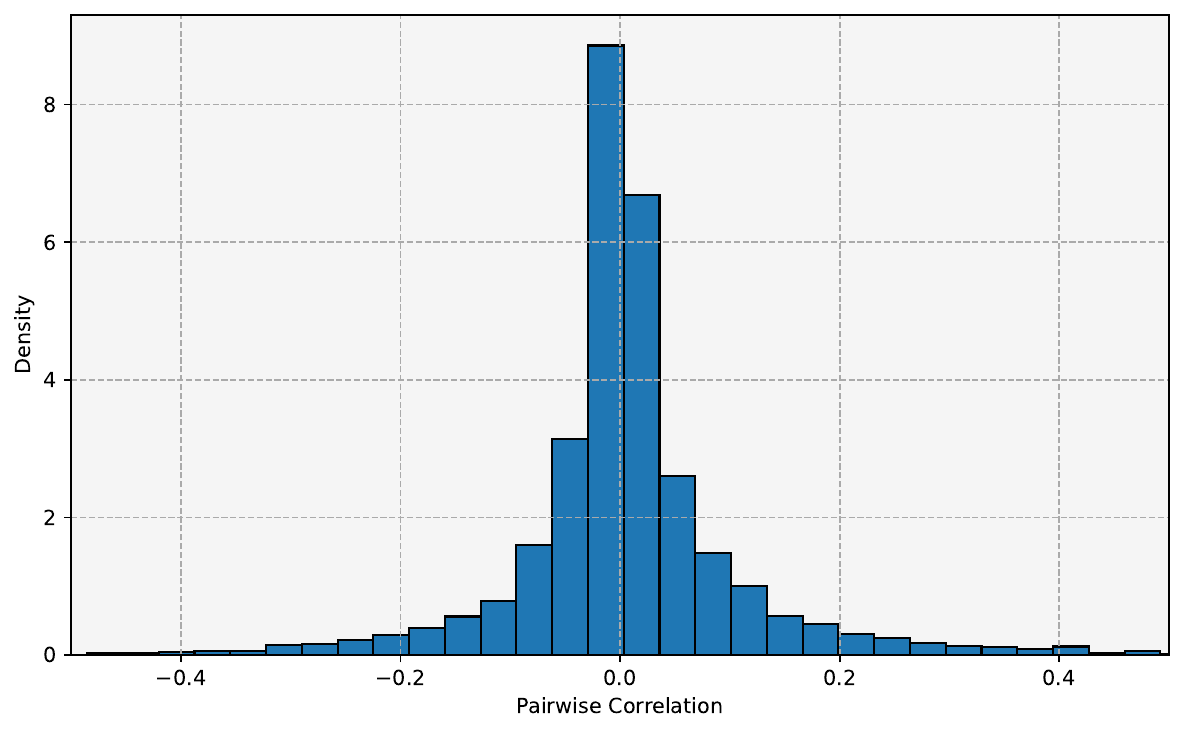}
    \caption{Density of pairwise correlations among the 102 firm characteristics.}
    \label{fig:corr_density}
\end{figure}

As a first step, we examine the correlation structure of the 102 firm characteristics used in the empirical analysis. 
\autoref{fig:corr_density} plots the empirical distribution of the 5,151 pairwise correlations, and \autoref{tab:corr_summary} reports the corresponding summary statistics. 
The distribution is centered close to zero, with mean $0.0079$ and median $0.0003$. 
Most pairwise correlations are small in magnitude, although moderate dependence remains in the tails. 
These summaries do not verify the strong irrepresentable condition, which depends on the full blockwise correlation structure rather than pairwise correlations alone, but they suggest that the empirical design is not dominated by pervasive near-collinearity.

\begin{table}[htbp]
\centering
\caption{Summary statistics of pairwise correlations among the 102 firm characteristics}
\label{tab:corr_summary}
\footnotesize
\setstretch{1.15}
\begin{tabular}{lccccccc}
\toprule
Statistic & Q5 & Q25 & Q50 & Q75 & Q95 & Mean & Std. \\
\midrule
Value & $-0.1563$ & $-0.0288$ & $0.0003$ & $0.0349$ & $0.1937$ & $0.0079$ & $0.1174$ \\
\bottomrule
\end{tabular}
\end{table}

\subsubsection{Simulation Diagnostic}

We next conduct a simulation diagnostic in the spirit of   \citet{zhao2006model}. 
The goal is to quantify how restrictive the strong irrepresentable condition is under dependence structures comparable to those in our data, and to distinguish between its global and entrywise implications.

Specifically, we generate $\bm X \in \mathbb{R}^p$ from a multivariate normal distribution with mean zero and correlation matrix $\bm C_{p\times p}$. 
To mimic the empirical dependence pattern in our sample, we construct $\bm C_{p\times p}=(c_{ij})$ by drawing the off-diagonal entries from the empirical distribution of pairwise correlations, setting the diagonal entries to one, and applying a positive-definite adjustment when needed. 

The first $q$ coordinates of $\bm X$ are treated as relevant and the remaining $p-q$ coordinates as inactive. 
We consider $p \in \{50,100,150,200\}$ and $q \in \{5,10,15,20\}$, with 100 replications for each $(p,q)$ pair. 
These values are chosen to reflect the empirical scale of the application, where the number of candidate characteristics is about $p=102$ and the number of selected characteristics is typically around $10$--$20$.

We report diagnostics for two designs. 
The block labeled \emph{Original Design} evaluates the strong irrepresentable condition for the raw characteristic vector $\bm X$. 
This serves as a baseline diagnostic for the original dependence structure. 
The block labeled \emph{Expanded Design} repeats the same exercise after transforming $\bm X$ into the augmented design matrix $\bm Z$ used in Step~1 of S-LASSO, including the linear terms, nonlinear basis terms, and interaction basis terms. 
This second diagnostic is directly aligned with our theory, since Assumption~3 is formulated for the expanded covariance matrix rather than for the raw design.

For each replication and each design matrix, let $\mathcal A$ and $\mathcal A^c$ denote the active and inactive component sets. 
Let
\begin{equation}
\bm r
=
\left|
\bm V_{\mathcal A^c\mathcal A}
\bm V_{\mathcal A\mathcal A}^{-1}
\mathrm{sign}(\bm\beta_{\mathcal A})
\right|,
\end{equation}
where the absolute value is taken componentwise. 
For an inactive component $j\in\mathcal A^c$, let $r_j$ denote the corresponding component of $\bm r$.
\begin{enumerate}
    \item The \emph{global satisfaction rate} is the fraction of replications in which $\max_{j\in\mathcal A^c} r_j<1$, so that the sample irrepresentable inequality holds for the entire inactive set. 
    \item The \emph{entrywise satisfaction rate} is the average fraction of inactive components satisfying $r_j<1$ across replications. 
\end{enumerate}
Thus, the global rate corresponds to the formal all-component requirement, whereas the entrywise rate is a diagnostic for whether violations are widespread or concentrated in a small number of inactive components.

\begin{table}[htbp]
\centering
\caption{Satisfaction rates of the strong irrepresentable condition under empirically calibrated dependence: diagnostics for the original design and the basis-expanded design used in S-LASSO.}
\label{tab:irrepresentable_simulation}
\footnotesize
\tabcolsep=0.25cm
\setstretch{1.5}
\begin{tabular}{llccccccccc}
\toprule
       &  & \multicolumn{4}{c}{\textbf{Original Design}} &  & \multicolumn{4}{c}{\textbf{Expanded Design}} \\
\cmidrule(lr){3-6} \cmidrule(lr){8-11}
       &  & $p=50$ & $p=100$ & $p=150$ & $p=200$ &  & $p=50$ & $p=100$ & $p=150$ & $p=200$ \\
\midrule
\rowcolor{gray!25} \multicolumn{11}{l}{\textbf{Panel A: Global Satisfaction Rate}} \\
$q=5$  &  & $94\%$    & $98\%$    & $100\%$   & $100\%$   &  & $94\%$    & $98\%$    & $100\%$   & $100\%$   \\
$q=10$ &  & $68\%$    & $81\%$    & $93\%$    & $95\%$    &  & $68\%$    & $81\%$    & $93\%$    & $95\%$    \\
$q=15$ &  & $27\%$    & $49\%$    & $61\%$    & $82\%$    &  & $27\%$    & $49\%$    & $61\%$    & $82\%$    \\
$q=20$ &  & $9\%$     & $20\%$    & $33\%$    & $46\%$    &  & $9\%$     & $20\%$    & $33\%$    & $46\%$    \\
\addlinespace[0.4ex]
\rowcolor{gray!25} \multicolumn{11}{l}{\textbf{Panel B: Entrywise Satisfaction Rate}} \\
$q=5$  &  & $99.84\%$ & $99.98\%$ & $100.00\%$ & $100.00\%$ &  & $99.96\%$ & $99.99\%$ & $100.00\%$ & $100.00\%$ \\
$q=10$ &  & $98.90\%$ & $99.73\%$ & $99.95\%$  & $99.97\%$  &  & $99.73\%$ & $99.93\%$ & $99.99\%$  & $99.99\%$  \\
$q=15$ &  & $95.57\%$ & $98.94\%$ & $99.59\%$  & $99.90\%$  &  & $98.89\%$ & $99.74\%$ & $99.90\%$  & $99.97\%$  \\
$q=20$ &  & $91.43\%$ & $97.46\%$ & $99.12\%$  & $99.54\%$  &  & $97.86\%$ & $99.37\%$ & $99.78\%$  & $99.89\%$  \\
\bottomrule
\end{tabular}
\end{table}

\autoref{tab:irrepresentable_simulation} reports two complementary measures. 
Panel~A gives the \emph{global satisfaction rate}, namely, the fraction of replications in which the strong irrepresentable condition holds for the entire inactive set. 
This is a worst-case requirement: a single violating inactive component is enough for the whole design to fail. 
Panel~B gives the \emph{entrywise satisfaction rate}, defined as the fraction of inactive components that satisfy the corresponding componentwise inequality. 
This is only a diagnostic summary and should not be interpreted as a substitute for the formal global condition.

Two findings are worth noting. 
First, the global satisfaction rate declines as $q$ increases, indicating that the condition becomes more demanding as the active set grows. 
Second, the entrywise satisfaction rate remains high in nearly all cases, even when the global rate is much lower. 
This suggests that global failure is often driven by a relatively small number of unfavorable inactive components rather than by widespread violations across the entire inactive set. 
Comparing the two blocks further shows that the basis expansion used in S-LASSO does not materially aggravate the condition in our empirically calibrated designs. 
In particular, the global diagnostic is essentially unchanged, while the entrywise diagnostic is often slightly improved after expansion. 
Hence, the main difficulty appears to come from the original dependence structure itself rather than from the basis expansion.


\end{document}